\documentclass[journal]{IEEEtran}
\ifCLASSINFOpdf
\else
\fi
\usepackage{amsmath} 
\usepackage{amssymb}  
\usepackage{graphicx}
\usepackage{amsmath}
\usepackage{amsfonts}
\usepackage{amssymb}
\usepackage{float}
\usepackage{color}
\usepackage{xcolor}
\usepackage{subcaption}
\usepackage[font=footnotesize]{caption}
\usepackage{balance}
\usepackage{cite}
\usepackage{multirow}
\usepackage{array}
\usepackage{booktabs,array,dcolumn}
\newcolumntype{P}[1]{>{\centering\arraybackslash}p{#1}}
\usepackage{mathtools}
\usepackage{algpseudocode}
\usepackage{listings}
\usepackage{cite}
\usepackage{color}
\usepackage{bm}
\usepackage{comment}
\usepackage[ruled,vlined]{algorithm2e}
\usepackage{amsmath}
\usepackage{graphicx}
\usepackage{color}
\usepackage{adjustbox}
\usepackage{centernot}
\usepackage{csquotes}
\usepackage{mathtools}
\usepackage{mathabx}
\usepackage{physics}
\usepackage{float}
\usepackage{multirow}
\usepackage{enumerate}
\usepackage{color}
\usepackage{cite}

\usepackage{amsmath} 
\usepackage{amssymb}  
\usepackage{bm}
\usepackage{xcolor}
\usepackage{hyperref}
\usepackage{etoolbox}
\newcommand{\rom}[1]{\uppercase\expandafter{\romannumeral #1\relax}}
\usepackage[normalem]{ulem}
\usepackage{balance}
\usepackage{bm}
\usepackage{balance}
\usepackage{mathtools}

\usepackage{amsthm}
\newtheorem{definition}{Definition}
\newtheorem{theorem}{Theorem}

\newtheorem{proposition}{Proposition}
\newtheorem{remark}{Remark}

\usepackage{setspace}

\usepackage{color}
\definecolor{gray}{RGB}{128,128,128}

\begin{document}

\title{ A Robust Stackelberg Game for Cyber-Security Investment in Networked Control Systems}
  \author{Pratishtha~Shukla,~\IEEEmembership{}
      Lu~An,~\IEEEmembership{}
      Aranya~Chakrabortty,~\IEEEmembership{Senior~Member,~IEEE,}
      and~Alexandra~Duel-Hallen,~\IEEEmembership{Fellow,~IEEE} 
\thanks{Manuscript received October 19, 2021; revised June 10, 2022 and July 30, 2022. This research was partially supported by the US National Science Foundation under grant ECCS 1544871.}
\thanks{Pratishtha Shukla was with the Department of Electrical and Computer Engineering, North Carolina State University, Raleigh, NC 27606 USA. She is now with the Oak Ridge National Lab, TN 37831 USA (e-mail: pshukla@ncsu.edu).}
\thanks{Lu An was with the Department of Electrical and Computer Engineering, North Carolina State University, Raleigh, NC 27606 USA. She is now with IBM Research, San Jose, CA 95120 USA (e-mail: lan4@ncsu.edu).}
\thanks{Aranya Chakrabortty and Alexandra Duel-Hallen are with the Department of Electrical and Computer Engineering, North Carolina State University, Raleigh, NC 27606 USA
(e-mail: achakra2@ncsu.edu; sasha@ncsu.edu).}}

\maketitle

\begin{abstract}
We present a resource-planning game for cyber-security of networked control systems (NCS). The NCS is assumed to be operating in closed-loop using a linear state-feedback $\mathcal{H}_2$-controller. A zero-sum, two-player Stackelberg game (SG) is developed between an attacker and a defender for this NCS. The attacker aims to disable communication  of selected nodes and thereby render the feedback gain matrix to be sparse, leading to degradation of closed-loop performance,  while the defender aims to prevent this loss by investing in the protection of targeted nodes. Both players trade their $\mathcal{H}_2$-performance objectives for the costs of their actions. The standard backward induction method is modified to determine a cost-based Stackelberg equilibrium (CBSE) that saves the players' costs without degrading the control performance. We analyze the dependency of a CBSE on the relative budgets of the players as well as on the node ``importance” order. Moreover, a robust-defense  method is developed for the realistic case when the defender is not informed about the attacker's resources. The proposed algorithms are validated using examples from wide-area control of electric power systems. It is demonstrated that reliable and robust defense is feasible unless the defender's resources are severely limited relative to the attacker's resources. We also show that the proposed methods are robust to time-varying model uncertainties and thus are suitable for long-term security investment in realistic NCSs. Finally, we employ computationally efficient genetic algorithms (GA) to compute the optimal strategies of the attacker and the defender in realistic large power systems.
\end{abstract}

\IEEEpeerreviewmaketitle

\section{Introduction}
Cyber-physical security of networked control systems (NCS) is a critical challenge for the modern society \cite{McLaughlin2016,Dibaji2019,Pasqualetti2013TAC,Mo2012,Book-DOS}.  While research on NCS security has focused on false data injection and intermittent denial-of-service (DoS) attacks \cite{Dibaji2019,Mo2012,Book-DOS,Li2018TAC, Basar2019CPS-Survey, ZhuBasar2015Mag, Agupta2021} or stealing information from the cyber system using advanced persistent threats (APTs) \cite{Yuan_new2020},
malicious destruction  of communication hardware (e.g. circuit boards, memory units, and communication ports) \cite{McLaughlin2016} or persistent distributed DoS (DDoS) attacks (where selected targets are flooded with messages so that they are unable to perform their services) \cite{DDOS_link} have received relatively little attention. In reality, these attacks can cause more severe damage to the communication network of an NCS compared to data tampering and intermittent jamming as they tend to disable communication and thus prevent feedback control for an extended period of time,  requiring expensive repairs or recovery efforts \cite{McLaughlin2016,DDOS_link}. 

A legitimate question, therefore, is how can network operators invest money for securing the important assets in an NCS against attacks that disable communication permanently under a limited budget? The same question applies to attackers in terms of targeting the best set of devices whose failure to communicate will maximize damage. These types of questions are best answered using game theory, which has been used as a common tool for modeling and analyzing cyber-security problems as it effectively captures conflicting goals of attackers and defenders \cite{Li2018TAC,Basar2019CPS-Survey, ZhuBasar2015Mag,Yuan_new2020, pirani2021game,yang2021modeling,Amin2013}. 
However, game-theoretic research for NCS security is often  unrelated to the model of the physical system and/or control methods \cite{Basar2019CPS-Survey,Li2022edgeAI,zhang2022riskmitigation}. These works also do not consider persistent DDoS or hardware attack models \cite{Yuan_new2020, Basar2019CPS-Survey, Li2022edgeAI} and tend to employ repeated games where the players update their investment strategies in response to the actions of their opponents in real time \cite{Basar2019CPS-Survey,zhang2022riskmitigation}. The latter approach, however, is not practical when a long-term, fixed security investment is required. Recently, a mixed-strategy (MS) investment game  for mitigation of hardware attacks on an NCS was presented in \cite{PSACC2019}. However, MS games \cite{ZhuBasar2015Mag, Balcan2015, Basar2019CPS-Survey} are also unsuitable for realistic, long-term security investment since they have randomized strategies and must be played many times to realize the expected payoffs. A long-term security investment game has been proposed in \cite{Lu2020}, but this game does not address NCS performance objectives. 

In this paper, we develop a Stackelberg game \cite{StackelbergvsCournot} for persistent malicious attacks on NCSs where fixed, non-randomized investment strategies are determined for both players. The NCS is assumed to be operating in closed-loop using a state-feedback $\mathcal{H}_2$-controller. The actions of the players are modeled as discrete investment levels into the  network nodes, which indicate the levels of effort and the resulting chances of success of attack and protection at each node.
The need for feedback control in the model guides the selection of the levels of attack and security investment at each node. 
The attacker aims to disable communication to/from a set of selected nodes, 
 which makes the feedback gain matrix sparse, thereby degrading the closed-loop $\mathcal{H}_2$-performance. The defender, on the other hand, invests in tamper-resistant devices \cite{Mo2012}, intrusion monitoring,  threat management systems that combine firewalls and anti-spam techniques \cite{Book-DOS}, devices or software that ensure authorized and authenticated access via increased surveillance \cite{ENSA2011}, etc. to prevent the attacks and maintain the  optimal $\mathcal{H}_2$-performance. A {\it Stackelberg Equilibrium }(SE) \cite{StackelbergvsCournot} of this game describes an optimal resource allocation of the two players given their respective budgets.   Moreover, the traditional {\it Backward Induction} algorithm is modified to compute a {\it cost-based} SE (CBSE), which saves the players' costs without compromising  their payoffs. We analyze the dependency of the players' payoffs at CBSE on the budgets and numbers of investment levels. In addition, to address the scalability issue of the proposed game for large networked systems, we enhance the bidirectional, parallel, evolutionary, genetic algorithm (BPEGA) \cite{Lu_thesis}, thus providing a computationally-efficient approach to finding a~CBSE.

Furthermore, 
the model parameters of an NCS usually vary over time due to changes in operating conditions, thereby making the NCS model  uncertain \cite{UncertainNCSBook}. Modifying the security investment as these changes occur can be time-consuming and expensive. Thus,  unlike in repeated dynamic games \cite{Basar2019CPS-Survey}, we seek {\it fixed}, {\it long-term} investment that provides  {\it robustness}  to model uncertainty.
Moreover, for a realistic setting  where the defender, who is the leader of the SG, does not know the follower's, i.e. the attacker's, capabilities, we develop a robust-defense  sequential algorithm  where the  defender protects against the most powerful hypothetical attacker.

  The proposed games are validated using  an example of wide-area oscillation damping  control for the IEEE 39-bus model, which represents the New England power system. First, we show that as the cost of defense per node increases, CBSEs of the proposed cost-based Stackelberg game (CBSG)  reveal the ``{\it important}" physical nodes
 \cite{Mihailo2013TAC,PSACC2019}, which are prioritized for protection and attack due to their impact on the control performance.  We also show that reliable control performance can be maintained unless the defender's resources are much more limited than the attacker's.  Second, we
demonstrate feasibility of robust, long-term protection for power systems with model uncertainty as well as the defender's uncertainty about the attacker's resources.  Third, we demonstrate the efficiency of the enhanced BPEGA method. Finally, to  show the applicability of the proposed techniques to large-scale networks, we extend our simulation study to the IEEE 68-bus system, which has more than twice the number of states than the 39-bus model.

For our problem, we consider a linear state feedback controller with a $\mathcal{H}_2$-control objective. The proposed security investment strategies are assumed to be fixed over long periods of time. The NCS model may change over such periods, but we show that our game is robust to model variations provided the change is estimated and followed by a new $\mathcal{H}_2$-control design.
This assumption, although idealistic, provides baseline performance bounds and allows for deployment of efficient computational tools developed for sparse $\mathcal{H}_2$ feedback controllers \cite{Mihailo2013TAC, Dorfler2014}. While more realistic constraints and stability guarantees were recently addressed for sparse control designs in, e.g., \cite{Lian2018sparsity,LoCicero2021Sparsity,Arastoo2016CDCsparsity,BAHAVARNIA2017sparsity}, these methods are computationally complex and do not easily extend to large NCSs. However, they can be readily adopted by our investment methods when more computationally efficient solutions are found.

{\it The main contributions of this paper are summarized below:}
\begin{itemize}
    \item Development of a Stackelberg security investment game for hardware or persistent DDoS attacks that deactivate communication and thereby degrade feedback control performance.
    \item Formulation and analysis of robust, long-term security investment methods for NCSs with uncertain dynamic models and  defender's uncertainty about the attacker's resources.
    \item Identification of critical assets of NCSs using sparse feedback control and proposed security methods. 
    \item  Utilization of genetic-algorithm-based numerical solutions to the proposed methods for applications to large-scale NCSs.
\end{itemize}

The rest of the paper is organized as follows. In Section II, we describe an NCS model as well as  attack and defense models.  The proposed  cost-based Stackelberg game for NCS security investment is presented  in 
Section III.  Section IV describes realistic, long-term, robust security investment methods. Numerical results are provided in Section V. Section VI concludes the paper.

\vspace{-0.2cm}
\section{Problem Formulation}
\subsection{Networked Control System Model and Controller}
\begin{figure}[t]
\centering
\includegraphics[width=\linewidth]{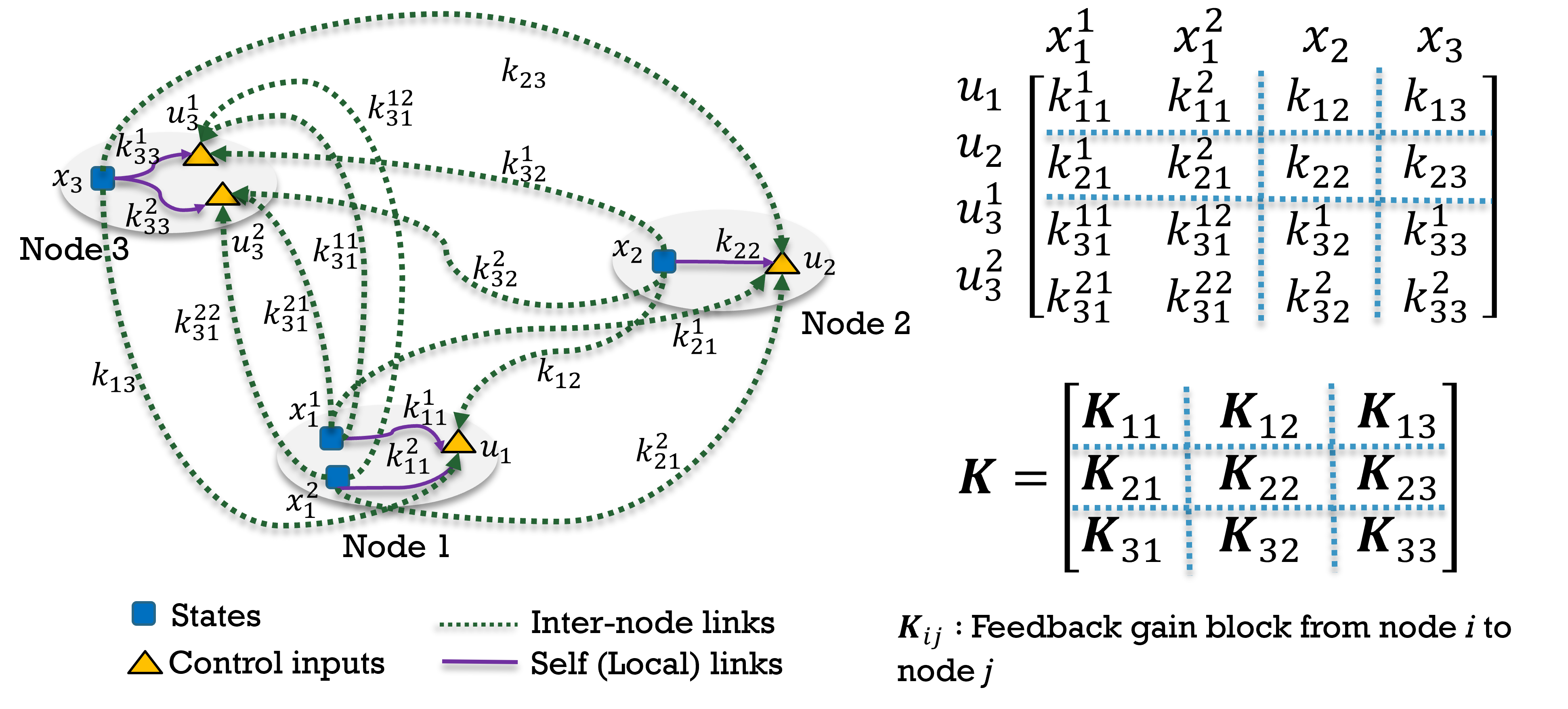}
\caption{ Schematic of a Networked Control System (NCS)}\label{fig1}
\end{figure}

We consider an NCS with $n$ nodes. Each node may be characterized by multiple states and control inputs, as shown in Fig. \ref{fig1}. 
At the $i^{th}$ node,  the state vector is denoted as $\bm{x}_i \in \mathbb R^{m_i}$,   with the total number of states $m = \sum_{i=1}^n m_i$, and the control input is denoted as $\bm{u}_i \in \mathbb R^{r_i}$, $i=1,\dots,\,n$, with the total number of control inputs $r = \sum_{i=1}^n r_i$. The state-space model of the network is written as
\begin{equation} \label{eq1}
\bm{\dot{x}}(t) = \bm{Ax}(t) + \bm{Bu}(t) + \bm{Dw}(t),
\end{equation}
where $\bm{x}(t)=(\bm{x}_1^T(t),...,\bm{x}_n^T(t))^T \in \mathbb{R}^{m\times 1}$, $\bm{u}(t)=(\bm{u}_1^T(t),...,\bm{u}_n^T(t))^T \in \mathbb{R}^{r\times 1}$, $\bm{w}(t) \in \mathbb{R}^{q\times 1}$ is a disturbance input modeled as white noise and $\bm{A} \in \mathbb{R}^{m \times m}$, $\bm{B} \in \mathbb{R}^{m \times r}$, $\bm{D} \in \mathbb{R}^{m \times q}$ are the state, input, and the disturbance matrices, respectively.  Assuming all states are measured, the control input $\bm{u}(t)$ is designed using  linear  state-feedback 
\begin{equation} \label{eq2}
\bm{u}(t)=-\bm{K}\bm{x}(t),
\end{equation}
where $\bm{K} \in \mathbb{R}^{r\times m}$ is the feedback gain matrix: 

\begin{equation}
\bm{K}= \begin{bmatrix}
\bm{K}_{11} & \bm{K}_{12} & \cdots & \bm{K}_{1m}\\
\bm{K}_{21} & \bm{K}_{22} & \cdots & \bm{K}_{2m}\\
\vdots & \vdots & \cdots & \vdots\\
\bm{K}_{r1} & \bm{K}_{r2} & \cdots & \bm{K}_{rm}\\
\end{bmatrix}.\label{eq3}
\end{equation}
From (\ref{eq3}) it follows that
\begin{equation} \label{eq4}
\bm{u}_j(t)=-\bm{K}_{ji}\bm{x}_i
\end{equation}
i.e. the $r_i \times m_j$ block matrix $\bm{K}_{ji}$ represents the topology of the communication network needed to transmit the state of node $i$ to the controller at node $j$. The diagonal blocks $\bm{K}_{ii}$ correspond to the local, or self-links, while the off-diagonal blocks $\bm{K}_{ji}$, $i \neq j$ indicate the inter-node communication links, respectively, as shown in Fig. \ref{fig1}.

The  objective is to find the feedback matrix $\bm{K}$ that minimizes the $\mathcal{H}_2$-performance cost function 
\begin{gather} \label{eq5}
 J(\bm{K}) = trace(\bm{D}^T\bm{P}\bm{D}) \\
\text{subject to} \ (\bm{A}-\bm{B}\bm{K})^T\bm{P}+ \bm{P}(\bm{A}-\bm{B}\bm{K}) \nonumber \\
\quad \quad \quad \quad \quad \quad \quad = -(\bm{Q}+\bm{K}^T\bm{R}\bm{K})\nonumber
\end{gather}
where $\bm{P}$ is the closed-loop observability Gramian, $\bm{Q}=\bm{Q}^T \succeq 0 \in \mathbb{R}^{m \times m}$ and $\bm{R}=\bm{R}^T \succ 0 \in \mathbb{R}^{r \times r}$ are design matrices that denote the state and control weights, respectively.  Using the standard assumptions, $(\bm{A},\bm{B})$ is stabilizable, and $(\bm{A},\bm{Q}^{1/2})$ is detectable
\cite{Mihailo2013TAC}.

\vspace{-0.15cm} 
\subsection{Attack and Defense Model} 
We consider two types of attacks, namely hardware attacks and persistent DDoS attacks. In hardware attacks, a malicious agent attempts to compromise the computer hardware components associated with selected network nodes, such as communication ports, circuit boards, etc. \cite{McLaughlin2016}, causing disruption of the feedback signal to/from these nodes. On the other hand, in persistent DDoS attacks, the computers associated with selected network nodes are flooded with messages and thus are unable to communicate control data for long time periods \cite{DDOS_link}. 
As shown in Fig. 1, the nodes in an NCS communicate with each other according to the connectivity of $\bm{K}$ to exchange the state information and compute the control inputs following (\ref{eq2}). The communication between the nodes of an NCS can be defined by various topologies, such as direct communication, internet or  an ad-hoc network, a cellular network, cloud-based or edge-based hierarchical communication, etc \cite{ge2017distributed}. 
Independent of the network model, the hardware components that enable transmission of each state $x_i(t)$ and reception of each state vector $u_i(t)$ from all network nodes (e.g., on-board computers, virtual machines in the cloud \cite{Soudbakhsh2014CDC}, etc) can be destroyed or compromised during a hardware attack \cite{McLaughlin2016}. Similarly, as discussed earlier, a persistent DDoS attack that targets node $i$ can prevent this node from sending and receiving feedback data for an extended time period.
\begin{figure}[t]
\centering
\includegraphics[width=\linewidth]{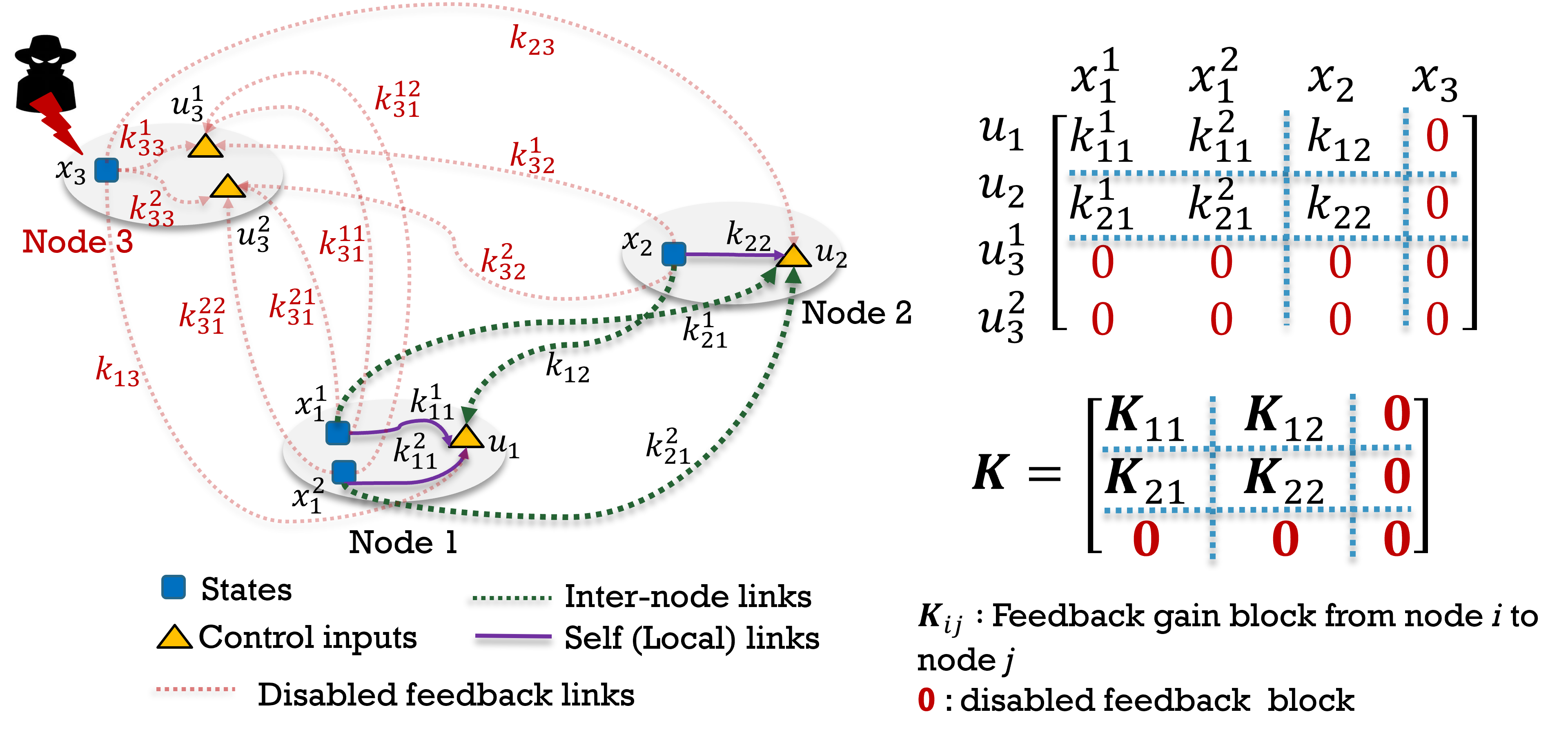}
\caption{Disabled communication for node 3 of Fig. \ref{fig1} corresponding to the sparsity pattern $s_6 =(110)$. All communication and self-feedback links associated with node 3 are disabled.}\label{fig2}
\end{figure}
In either case, a successful attack on node $i$ will zero out the entire $i^{th}$ block-row and block-column of $\bm{K}$ in (\ref{eq3}) as shown in Fig. \ref{fig2}, thereby degrading
the closed-loop $\mathcal{H}_2$-performance (i.e. resulting in a suboptimal value of $J$ in (\ref{eq5})). 
The attacker invests as per its budget into selected nodes to increase the value of $J$  while the defender aims to protect  the system performance from such degradation against attacks by installing tamper-resistant devices and/or intrusion monitoring software (for details of the mechanisms please see \cite{Mo2012,Book-DOS,ENSA2011}). Note that the defender does not know  when and where an attack might happen, so it acts {\it proactively} by selecting a set of nodes and the protection levels to maintain $J$ as low as possible within the defense budget in case of a future attack. Our objective is to formulate a Stackelberg game for optimal resource allocation by both the defender of an NCS and a malicious attacker assuming the hardware and persistent DDoS attack models described~above.

\vspace{-0.1cm}
\section{ Stackelberg Game for NCS Security Investment}
The attacker's actions $\bm{a}$ and the defender's actions $\bm{d}$ indicate the levels of investment into the system nodes, which measure the levels of effort and the resulting chances of successful attack and protection at these nodes, respectively. Both players have budget constraints and also aim to reduce their costs of attack or protection. The attacker and the defender have opposite control performance goals. The attacker aims to increase the $\mathcal{H}_2$-performance cost in (\ref{eq5}) while the defender tries to keep it as close to the optimal value as possible. The payoffs of the attacker and defender are denoted as $U^a(\bm{a},\bm{d})$ and $U^d(\bm{a},\bm{d})$, respectively. Given a pair of investment strategies $(\bm{a},\bm{d})$, $U^d(\bm{a},\bm{d}) = -U^a(\bm{a},\bm{d})$, resulting in a {\it zero-sum} game~\cite{Osborne1994}. 
 
As in many SGs for security \cite{Li2018TAC,Lu2020, Yuan_new2020}, the defender is the {\it leader} and chooses its  investment profile  first. Given a defender's action $\bm{d}$, the attacker follows by a best-response to $\bm{d}$, given by  $\bm{a} = g(\bm{d}) =\mathop {\arg \max} \limits_{\bm{a}} U^a(\bm{a},\bm{d})$. Thus, the defender chooses a strategy $\bm{d}$ that maximizes its payoff given the attacker's best responses to its actions. A resulting  {\it Stackelberg Equilibrium} (SE) \cite{StackelbergvsCournot} specifies a pair of strategies $(\bm{a}^*,\bm{d}^*)$, which optimizes the payoffs of the players in an SG.
  Finally, we augment the standard SG described above by selecting an SE that reduces the players' costs. 
  The resulting game is termed  as {\it cost-based Stackelberg game} (CBSG).
  
In this section, we assume that the opponent's budget and the number of investment levels are known to each player. Moreover, we assume that the system model is fixed and known to the players. These idealistic assumptions result in a baseline game performance characterization and will be relaxed in Section IV.
 \subsection{Player's Actions and Cost Constraints}
The actions of the players are given by {\it n-}dimensional
investment vectors, denoted as 
\begin{equation}\label{eq6}
\bm{a} = (a_1 , a_2, ...,a_n) , \ \bm{d}= (d_1,d_2,...,d_n) 
\end{equation}
for the attacker and the defender, respectively. A higher value of $a_k$ (or $d_k$) corresponds to a larger attack (or protection) investment level, thus the level of effort, at node $k$. For example, increased effort $a_k$ exerted by the attacker over the node $k$ leads to increased loss (e.g., through increased probability of compromising this node) while the defender's action $d_k$ over the same node is to choose a level of protection, effort, or investment in security, in order to mitigate the attack. The greater the defender's investment, the lower the effectiveness of the attacker's effort. Investing at level $d_k = 1$ is assumed to provide node $k$ with perfect protection \cite{Sarabi2014}. The levels $a_k$ in (\ref{eq6}) are chosen from the set  $\left\{0,\frac{1}{L_a},\frac{2}{L_a},...,1\right\}$ where $L_a+1$ is the total number of attacker's investment levels. Similarly, $d_k \in \left\{0,\frac{1}{L_d},\frac{2}{L_d},...,1\right\}$ where $L_d + 1$ is the number of defender's investment levels. Given the actions $\bm{a}$ and $\bm{d}$, the probability of successful attack at node $k$ is given by \cite{Sarabi2014}
\begin{equation}\label{eq8}
P_k (\bm{a},\bm{d})= a_k(1-d_k).
\end{equation}

The set of possible attack outcomes at all nodes is represented by a set of $2^n$ {\it sparsity patterns}, or binary {\it n-}tuples, 
\begin{equation}\label{eq9}
\bm{s}_m = (s_m^1,...,s_m^n) 
\end{equation}
where $s_m^k = 0$ indicates an attack is successful at node $k$ while $s_m^k = 1$ means that either protection is successful at node $k$ or node $k$ is not attacked. From (\ref{eq8}), the probability that the sparsity pattern $\bm{s}_m$ occurs given the strategy pair $(\bm{a},\bm{d})$ is 
 \begin{equation}
       \label{eq10}
      P_{\bm{s}_m}(\bm{a},\bm{d})  = \prod_{k, \ s_m^k = 0}P_k(\bm{a},\bm{d}) \ \prod_{k, \ s_m^k = 1}(1-P_k(\bm{a},\bm{d})).
\end{equation}
   
Finally, the players' cost constraints are as follows. The bounds on the attacker's and defender's budgets are denoted by $R_a$ and $R_d$, respectively. Let $g_{a_i}$ and $g_{d_i}$ be the costs of attacking and protecting node $i$ at full effort, respectively. Scaling this cost by the level of effort and summing over all nodes, the actions of the players are cost-constrained as 
\begin{equation}\label{eq11a}
         \sum_{i=1}^n g_{a_i}a_i  \leq R_a,         \sum_{i=1}^n g_{d_i}d_i  \leq R_d. 
\end{equation}
We normalize (\ref{eq11a}) by dividing the first and the second inequality by $R_a$, $R_d$, respectively, with the normalized cost per node at full effort given by $\gamma_{a_i} = \frac{g_{a_i}}{R_a}$ and $\gamma_{d_i} = \frac{g_{d_i}}{R_d}$.
Thus, the normalized cost constraints are given by
\begin{equation} \label{eq11aa}
    \sum_{i=1}^n \gamma_{a_i}a_i  \leq 1,  \quad \sum_{i=1}^n \gamma_{d_i}d_i  \leq 1 
\end{equation}

\subsection{Structural Sparsity and Players' Payoffs}
Following an attack, when a sparsity pattern $\bm{s}_m$ as in (\ref{eq9}) occurs, all communication to/from each node $k$ for which $s_m^k = 0$ is disabled. Thus, the corresponding feedback matrix $\bm{K}$ in (\ref{eq3}) has the sub-blocks $\bm{K}_{kp} = \bm{0}$ and $\bm{K}_{qk} = \bm{0}$ for all $p = 1,2,...,n, \ q = 1,2,...,n$, imposing the {\it structural sparsity} constraint \cite{Mihailo2013TAC} on the matrix $\bm{K}$. For example, Fig. \ref{fig2} shows the scenario where communication is disabled within and to/from node 3 and the resulting structural sparsity of the feedback matrix $\bm{K}$.
 
Next, we define the {\it $\mathcal{H}_2$-performance loss} vector $\bm{\Delta} = (\Delta_{\bm{s}_0},..,\Delta_{\bm{s}_{m}},...,\Delta_{\bm{s}_{2^n-1}})$ with the $m^{th}$ element given by 
\begin{equation}\label{eq12}
\Delta_{\bm{s}_{m}} = J(\bm{K}^*_{\bm{s}_{m}}) - J(\bm{K}^*_{\mathcal{H}_2}), 
\end{equation}
where $\bm{K}^*_{\mathcal{H}_2}$ and $\bm{K}^*_{\bm{s}_{m}}$ optimize  the $\mathcal{H}_2$-performance objective function in (\ref{eq5}) without and with the structural sparsity constraint imposed by $\bm{s}_{m}$, respectively. The latter can be computed using the structured $\mathcal{H}_2$ optimization algorithm in \cite{Mihailo2013TAC}.  The loss corresponding to the sparsity pattern $\bm{s}_0=(0,...,0)$ represents the open-loop loss given by
\begin{equation}\label{eq12l}
    \Delta_{OL} = J_{OL} - J(\bm{K}^*_{\mathcal{H}_2})
\end{equation}
 where $J_{OL}$ is $\mathcal{H}_2$-performance of the open-loop system.

Given $\bm{a}$ and $\bm{d}$, the attacker's payoff is
\begin{equation}
    \label{eq13a}
    U^a(\bm{a},\bm{d}) = \sum_{m=0}^{2^n-1} \ P_{\bm{s}_m}(\bm{a},\bm{d}) \Delta_{\bm{s}_m}
\end{equation}
and the defender's payoff is
\begin{equation}\label{eq13b}
U^d(\bm{a},\bm{d}) = - U^a(\bm{a},\bm{d}),
\end{equation}
where $P_{\bm{s}_m}(\bm{a},\bm{d})$ is given by (\ref{eq10}). 
Thus, the attacker and the defender aim to maximize and minimize, respectively, the expected system loss in (\ref{eq13a}) under the constraints (\ref{eq11aa}). 
Given the two payoffs, we next define the Cost-based Stackelberg Equilibrium of our proposed game.
\begin{remark}
 The proposed investment methods can easily adopt a different controller type or control metric by substituting another control objective instead of the $\mathcal{H}_2$-metric {\normalfont (\ref{eq5})} into the loss expression {\normalfont (\ref{eq12})} (and {\normalfont (\ref{eq12l})}). While we assume stable open loop and thus stabilizing optimal sparse feedback controllers as in the wide-area control  example in Section V, the game can be extended 
 to scenarios where the attack leads to instability by protecting the corresponding sparsity patterns 
fully.
\end{remark}


\subsection{Cost-Based Stackelberg Equilibrium (CBSE)}
Several methods can be employed for computing an SE \cite{simaan1973}. Among these, the {\it Backward Induction} (BI) algorithm is a popular approach for finite SGs \cite{Li2018TAC,Yuan_new2020}. Since multiple SEs are possible in the proposed SG, we modify the BI method to select an SE that saves both players' investment costs \cite{Lu2020,Lu_thesis}. The resulting {\it Cost-Based Backward Induction} (CBBI) method is summarized in Algorithm \ref{alg:ALG1} where the steps 1(a) and 2(b) are included to save the attacker's and the defender's costs without reducing their payoffs. 

\begin{algorithm} [ht]
    \SetAlgoLined
            \caption{Cost-Based Backward Induction Algorithm} 
            \label{alg:ALG1}
{\bf Step 1:} For each action $\bm{d}$ by the defender that satisfies (\ref{eq11aa})\;

{\bf (a)} {\it the attacker chooses its best response $g(\bm{d})$} under its cost constraint in (\ref{eq11aa}):
\begin{align}\label{eq14}
    g(\bm{d}) = \mathop {\arg \max }\limits_{\bm{a}} U^a(\bm{a},\bm{d}), \quad
    \mbox{s.t} \ \sum_{i=1}^n \gamma_{a_i}a_i \leq 1. 
    \end{align}

{\bf (b)} If there are multiple best responses $g(\bm{d})$ in Step 1(a), the attacker chooses the best response with the {\it smallest cost}:
\begin{equation}
    \label{eq15}
    g_o(\bm{d}) = \mathop {\arg \min }\limits_{g(\bm{d})} (\sum_{i=1}^n \gamma_{a_i}g(\bm{d})_i) .
\end{equation}

{\bf Step 2:} For defender-\;

{\bf (a)} Assuming the attacker uses a response (\ref{eq15}) for each defender's action $\bm{d}$, the {\it defender chooses an investment strategy that maximizes its payoff} under the constraint in (\ref{eq11aa}):
\begin{align} \label{eq16}
    \bm{d}^* = \mathop {\arg \max }\limits_{\bm{d}} U^d(g_o(\bm{d}),\bm{d}), \ \
\mbox{s.t} \ \sum_{i=1}^n \gamma_{d_i}d_i \leq 1.
\end{align}

{\bf (b)} If (\ref{eq16}) has multiple solutions, then a defender's strategy with the {\it smallest cost} is chosen:
\begin{equation}
    \label{eq17}
    \bm{d}^*_o = \mathop {\arg \min }\limits_{\bm{d}^*} (\sum_{i=1}^n \gamma_{d_i}d^*_i).
\end{equation}
Denote
\begin{equation}
    \label{eq18}
    \bm{a}_o^* = g_o(\bm{d}_o^*).
\end{equation}

\end{algorithm}

\begin{definition}
The strategy pair $(\bm{a}_o^*,\bm{d}^*_o)$ in {\normalfont (\ref{eq17})} and {\normalfont (\ref{eq18})} is a {\it Cost-based Stackelberg Equilibrium} (CBSE) of the proposed cost-based Stackelberg game (CBSG). 
\end{definition} 
\begin{remark}
Note that in Steps 1(b) and 2(b), ties are resolved arbitrarily. Moreover, at CBSE, $U^a(\bm{a}_o^*,\bm{d}_o^*) \geq 0$ since a CBSE cannot have a lower attacker's payoff than in the ``no attack" case when $U^a(\bm{0},\bm{d}) = 0$. Since the game is zero-sum, $U^d(\bm{a}_o^*,\bm{d}_o^*) \leq 0$. Finally, while in general SGs an SE might not exist, the CBBI Algorithm \ref{alg:ALG1} guarantees existence of a CBSE in the proposed finite CBSG (see Theorem 1(a)).
\end{remark}
\begin{remark}
 The complexity of computing $m^{th}$ element of the loss vector {\normalfont (\ref{eq12})} using structural optimization {\normalfont \cite{Mihailo2013TAC}} is $O(2^n \times nnz(\bm{K}_{\bm{s}_m}))$ where $nnz(.)$ represents the number of non-zero elements in the structured matrix $\bm{K}_{\bm{s}_m}$. In the rest of the paper, we refer to ``complexity" as the computational complexity of the actual game or algorithm steps after the loss vector is computed. 
\end{remark}
Since the number of possible payoffs of each player {\normalfont (\ref{eq13a})},{\normalfont (\ref{eq13b})} of the proposed CBSG is bounded by $(L_a +1)^n \times (L_d+1)^n$, the worst-case game complexity is $O((L_a +1)^n \times (L_d+1)^n)$. The actual number of payoffs can be significantly lower due to the cost constraints in {\normalfont (\ref{eq11aa})}, which limit the sets of player's actions. 

Following \cite{Lu2020}, we summarize some properties of the proposed game (\ref{eq13a})-(\ref{eq13b}) in Theorem~1. 
\begin{theorem} 
($\bf{a}$) Since the number of possible actions is finite, Algorithm \ref{alg:ALG1} has at least one solution (a CBSE). \newline
($\bf{b}$) The control performance loss {\normalfont (\ref{eq12})}  is the same for all solutions of the game (all CBSEs). \newline
($\bf{c}$) Given $L_a$ and $L_d$, the payoff of the attacker {\normalfont (\ref{eq13a})} or defender {\normalfont (\ref{eq13b})}  does not increase with its cost per node when the opponent's cost per node {\normalfont (\ref{eq11aa})} is fixed. \newline
($\bf{d}$) Given $L_a$ and $L_d$, there exist $\epsilon > 0$ and $\theta > 0$ such that when $\gamma_a < \epsilon$ while $\gamma_d > \theta$, the attacker's payoff at CBSE $U^a(\bm{a}^*,\bm{d}^*) = \Delta_{OL}$ (the open-loop loss {\normalfont (\ref{eq12l})}). Moreover, there exists an $\alpha > 0$ such that when $\gamma_d < \alpha$, the attacker's payoff at CBSE $U^a(\bm{a}^*,\bm{d}^*) = 0$ (i.e. the optimal $\mathcal{H}_2$-performance is achieved). \newline
($\bf{e}$) When $L_d+1$ (or $L_a+1$) is increased to $L_d'$ (or $L_a'$) that satisfies $L_d' = \eta(L_d )$ (or $L_a' =\eta(L_a)$), where $\eta$ is a positive integer, the defender's {\normalfont (\ref{eq13b})} (or attacker's {\normalfont (\ref{eq13a})}) payoff does not decrease if the costs per node of both players and the opponent's number of investment levels $L_a+1$ (or $L_d+1$) are fixed.
\end{theorem}
\begin{proof}
Please refer to Appendix A.
\end{proof}
From Theorem 1, CBBI finds an SE with reduced costs of both players while providing the payoff of any other SE. The properties of the CBSG summarized in Theorem 1 will be illustrated in Section V. We next present a numerical approach for computing a CBSE.

\subsection{Cost-based Bidirectional Evolutionary Method for Computing a CBSE}
 The traditional BI algorithm and thus the CBBI method (\ref{eq14})-(\ref{eq18}) in Algorithm \ref{alg:ALG1} are referred to as  traversal searching methods. According to Remark 3, the CBBI Algorithm \ref{alg:ALG1} has $O((L_a +1)^n (L_d +1)^n)$  computational complexity, which rapidly increases as the system size ($n$) and the players’ numbers of investment levels ($L_a$ and $L_d$) grow. To reduce the computational complexity of finding an SE, genetic algorithms (GA) have been developed for SGs in \cite{LIU1998,Liu2009}. However, in these papers, the evolutionary process was employed only by the defender. To overcome this limitation, we proposed~a bidirectional parallel evolutionary GA (BPEGA) method in~\cite{Lu_thesis}. 

\begin{algorithm}[t!]
\SetAlgoLined
\KwResult{One strategy pair $(\bm{a}^*,\bm{d}^*)$}
{\bf Parameter initialization:} Population sizes $S_a$ and $S_d$, Crossover probability $P_c$, Mutation rate $P_m$, Maximum number of generations $T$; Current generation $t = 0$\\ 
{\bf Step 1. Population Initialization:} Randomly selected feasible initial population for both players $POP^0_a = \left\{\bm{a}_1,\dots,\bm{a}_{S_a}\right\}$ and $POP^0_d = \left\{\bm{d}_1,\dots,\bm{d}_{S_d}\right\}$, where $\forall \bm{a} \in POP^0_a$  and $\forall \bm{d} \in POP^0_d$  satisfy (\ref{eq11aa}) \;
\While{the termination criteria are not satisfied,}{
{\bf Step 2. Evaluation:} The defender and attacker compute $U^d(\bm{a},\bm{d})$ and $U^a(\bm{a},\bm{d})$ for all $\bm{a} \in POP^0_a$  and $\forall \bm{d} \in POP^0_d$ and evaluate them to compute their fitness values-
\begin{equation}
    f_d(\bm{d}) = U^d(g_t(\bm{d}),\bm{d}), \ \forall \bm{d} \in POP^t_d \label{eq21fitd}
 \end{equation}
 where a best response $ g_t(\bm{d}) = \mathop {\arg \max }\limits_{\bm{a}\in POP_a^t} U^a(\bm{a},\bm{d})$ and 
  \begin{equation}
     f_a(\bm{a}) = U^a(\bm{a},\bm{d}'), \ \forall \bm{a} \in POP^t_a  \label{eq22fita}
 \end{equation}
 where $\bm{d}' = \mathop {\arg \max }\limits_{\bm{d}\in POP_d^t} f_d(\bm{d})$\;
 
\For{Attacker and Defender} {
{\bf Step 3. Selection:} Select $S_a/2$ (or $S_d/2$) pairs of parents $tmp_p^a$ (or $tmp_p^d$) using the Roulette Wheel selection method \cite{LIU1998} \\
{\bf Step 4. Reproduction:} Apply crossover with $P_c$ and mutation operation with $P_m$ to generate $S_a$ (or $S_d$) children $tmp_c^a$ (or $tmp_c^d$) \cite{WhitleyGA1994}\\
{\bf Step 5. Check feasibility:} For each individual in $tmp_c^a$ (or $tmp_c^d$), check if it is a feasible solution to (\ref{eq11aa}). Include all feasible children in the set $tmp_{c,f}^a$ (or $tmp_{c,f}^d$)\\
{\bf Step 6. Combine, Sort and Reorder:} Combine $POP_a^t$ (or $POP_d^t$) with the set of feasible children $tmp_{c,f}^a$ (or $tmp_{c,f}^d$), sort by the fitness value in the {\it descending order}. For all individuals with same fitness value, re-order them in ascending order of their costs using (\ref{eq11aa}). Individuals with same fitness value and cost are sorted according to their {\it seniority} of being selected, with the individuals who were present in an earlier generation placed ahead of those selected later. $S_a$ (or $S_d$) individuals with the highest rank are selected as $POP_a^{t+1}$ (or $POP^{t+1}_d$)\\
}
$t \leftarrow t+1$
}
{\bf Step 7.} Apply the CBBI Algorithm \ref{alg:ALG1} (\ref{eq14})-(\ref{eq18}) to the final generation $POP_a^T$ and $POP_d^T$ to determine $(\bm{a}^*,\bm{d}^*)$.
\caption{Cost-Based Bidirectional Parallel Evolutionary Genetic Algorithm}
\label{alg:ALG2}
\end{algorithm}

In this paper, we extend the BPEGA algorithm to find a CBSE of the proposed CBSG. The resulting method is termed cost-based BPEGA  and  is summarized in Algorithm \ref{alg:ALG2}. The changes from  the BPEGA Algorithm 5.2 of  \cite{Lu_thesis} are  as follows. 
 First, the elitist step (step 6) is modified to {\it reorder} the individuals with the same fitness value, which guarantees the selection of  the least-cost SE strategy. Furthermore, once the termination criteria described in Section 5.4.3 in \cite{Lu_thesis} are met, the CBBI Algorithm \ref{alg:ALG1} is applied on the final generation of populations of the players  in step 7 to find an optimal strategy pair $(\bm{a}^*,\bm{d}^*)$. Convergence of Algorithm \ref{alg:ALG2} follows from Proposition 5.1 in  \cite{Lu_thesis}.

\begin{proposition}
For the proposed CBSG, assume the crossover probability $P_c > 0$ and a mutation rate $P_m > 0$. The defender and attacker's payoffs generated by the CB-BPEGA (Algorithm \ref{alg:ALG2}) converge to their respective payoffs at a CBSE as the number of iterations tend to infinity.
\end{proposition}
\begin{proof}
 Please refer to Appendix B.
\end{proof}


Finally, we compare the computational complexity of the traversal CBBI Algorithm \ref{alg:ALG1} (\ref{eq14})-(\ref{eq18}) and the proposed CB-BPEGA. Given a fixed maximum number of generations $T$ and the population sizes $S_a$ and $S_d$ of the attacker and defender, respectively,  Algorithm \ref{alg:ALG2} has $O((T+1)S_aS_d)$ computational complexity (Section 5.4.4 in \cite{Lu_thesis}), which is significantly lower than the computational complexity of CBBI Algorithm \ref{alg:ALG1} (see Remark 3) when the system size is large and $TS_aS_d \ll (L_a+1)^n(L_d+1)^n$ \cite{Lu_thesis}.
In practice, NCSs with different system sizes might require different $T$ values to achieve convergence as demonstrated in Section V.

\section{Long-term, Robust Security Investment Methods for NCSs}
\subsection{Nominal-model CBSG for Systems with Model Uncertainty}
The CBSG developed above assumes a fixed system model in (\ref{eq1}), referred to below as the {\it nominal} model $M_{nom}$. In practice, the model varies with time. We refer to upcoming models as ``uncertain models" and denote the $i^{th}$ uncertain model as  $ M_{i}=\{\bm{A}_i,\,\bm{B}_i\}$ where $\bm{A}_i$ and $\bm{B}_i$ are the state and input matrices (\ref{eq1})  of this model. We set $M_1 = M_{nom}$.

If a CBSE $(\bm{a}^*_{nom},\bm{d}^*_{nom})$ of the SG designed for $M_{nom}$ is employed as an investment strategy  for a future system realization $M_i$, we refer to this game as the ``{\it nominal-model}" SG. Note that the
 CBSE of the latter game  is in general a suboptimal investment strategy for model $M_i$. The fractional difference
 between the payoffs of the model $M_i$ and the nominal model $M_1 = M_{nom}$ at the CBSE $(\bm{a}^*_{nom},\bm{d}^*_{nom})$ is defined as 
\begin{equation}\label{eq19}
\small
   \mu_{i,nom} \ \%= \abs{\frac{\bm{U}^a_{M_{i}}(\bm{a}^*_{nom},\bm{d}_{nom}^*) - \bm{U}^a_{M_{nom}}(\bm{a}^*_{nom},\bm{d}_{nom}^*)}{\bm{U}^a_{M_{nom}}(\bm{a}^*_{nom},\bm{d}_{nom}^*)}} \times 100 \ \%
\end{equation}
where $U^a_{M_i}(\bm{a},\bm{d})$ is the attacker's payoff (\ref{eq13a}) when the actual model is $M_i$.  This payoff is obtained as follows. First, define the loss vector $\bm{\Delta}^{M_i}$, which is found as in (\ref{eq12}) with $\bm{K}^*_{\bm{s}_m}$ and $\bm{K}^*_{\mathcal{H}_2}$ computed for model $M_i$. Then $U_{M_i}^a(\bm{a},\bm{d})$ and $U^d_{M_i}(\bm{a},\bm{d})$ are given by (\ref{eq13a}), (\ref{eq13b}), respectively, with $\bm{\Delta}$ replaced by $\bm{\Delta}^{M_i}$.
Note that $\mu_{1,nom} =0$.  A similar game can be designed based on an initial model $M_i$, $i\neq 1$, instead of $M_{nom}$.

\subsection{Robust-Defense Security Method}
In the above game, we made an idealistic assumption that the attacker and defender have complete information about each other.  In this subsection, we consider a realistic scenario where the defender, who acts first, does not know the attacker's resources and thus prepares to protect against the most powerful attacker. On the other hand, the attacker observes the defender's action and is thus able to generate its best response regardless of its knowledge about the defender's resources. The {\it robust-defense} method presented below can be formulated as a trivial {\it Bayesian Stackelberg game} (BSG) \cite{dynGame, PS_thesis}, but it is more intuitive to describe it as the following {\it sequential algorithm.} First, we assume the nominal model $M_1 = M_{nom}$ and drop the model index as in Section III. 
The utility of the defender is independent of $\bm{a}$ and is defined as:
\begin{equation}\label{eq25rdsg}
   U^d(\bm{d}) = -U^a(\bm{1}_n,\bm{d}),
\end{equation}
where $U^a(\bm{a},\bm{d})$ is given by (\ref{eq13a}), and $\bm{1}_n = (1,\cdots,1)_n$, i.e. the defender assumes full attack on all nodes. Once the defender identifies its optimal action $d^*_{RD}$ that maximizes (\ref{eq25rdsg}) while saving its cost, the attacker finds its cost-based best response using its utility (\ref{eq13a}).  The proposed  robust-defense (RD) method shown in Algorithm \ref{alg:ALG3} does not require backward induction since the players' optimal strategies are found sequentially.
\begin{algorithm} [htp]
    \SetAlgoLined
            \caption{Cost-Based Robust Defense Method } 
            \label{alg:ALG3}
{\bf A. Leader Optimization Subproblem:}

{\it Step 1(a)}- The defender chooses an investment by solving:
\begin{align}\label{eq:ch4.32}
    \bm{d}^* = \mathop {\arg \max }\limits_{\bm{d}} U^d(\bm{1}_n,\bm{d}), \quad 
\mbox{s.t} \ \sum_{i=1}^n \gamma_{d_i}d_i \leq 1.
\end{align}

{\it 1(b)}- A smallest cost defender's strategy is then chosen:
\begin{align} \label{ch4.33}
    \bm{d}^*_{RD} = \mathop {\arg \min }\limits_{\bm{d}^*} (\sum_{i=1}^n \gamma_{d_i}d^*_i).
\end{align}

{\bf B. Follower Optimization Subproblem:}

{\it Step 1(a)}- For given defender's strategy $\bm{d}^*_{RD}$, the attacker chooses its best response $\bm{a}^*$ solving:
\begin{align} \label{eq:ch4.34}
    \bm{a}^* = \mathop {\arg \max }\limits_{\bm{a}} U^a(\bm{a},\bm{d}_{RD}^*), \quad
    \mbox{s.t} \ \sum_{i=1}^n \gamma_{a_i}a_i \leq 1. 
\end{align}

{\it 1(b)}- A smallest cost attacker's strategy is chosen:
\begin{align} \label{eq:ch4.35}
    \bm{a}^*_{RD} = \mathop {\arg \min }\limits_{\bm{a}^*} (\sum_{i=1}^n \gamma_{a_i}a^*_i).
\end{align}
            
\end{algorithm}
 We denote a solution computed by Algorithm \ref{alg:ALG3} as $(\bm{a}^*_{RD},\bm{d}^*_{RD})$. The existence of this solution is guaranteed due to the finite action spaces of both players.

The GA implementation of the RD algorithm is a modified version of Algorithm \ref{alg:ALG2}, which was presented in  Algorithm 4.4 in \cite{PS_thesis}.   
Using a similar argument to the proof of Proposition 1, it can be easily shown that the payoffs of this sequential GA converge to the payoffs of a solution of the RD method as the number of iterations tends to infinity.

The attacker's payoff at the output of Algorithm \ref{alg:ALG3} is $U^a(\bm{a}^*_{RD},\bm{d}^*_{RD})$  (\ref{eq13a}). However, the actual payoff of the defender is not given by its optimized utility (\ref{eq25rdsg}) at $\bm{d}^*_{RD}$  since the actual attacker's action is not $\bm{a} = \bm{1}$, but is $\bm{a}^*_{RD}$. Thus, the actual defender's payoff is  
\begin{equation}\label{eqactualpayoffrdsg}
    U^d(\bm{a}^*_{RD},\bm{d}^*_{RD})=-U^a(\bm{a}^*_{RD},\bm{d}^*_{RD}),
\end{equation}
which is suboptimal due to the defender's uncertainty
as shown in Theorem 2.

\begin{theorem}
The defender's actual payoff $U^d(\bm{a}^*_{RD},\bm{d}^*_{RD})$ of the RD method {\normalfont (\ref{eqactualpayoffrdsg})} does not exceed its payoff at a CBSE of the CBSG {\normalfont (\ref{eq13a})-(\ref{eq13b})}, i.e. 
\begin{equation}
    U^d(\bm{a}^*_o,\bm{d}^*_o) \geq  U^d(\bm{a}^*_{RD},\bm{d}^*_{RD})
\end{equation}
while the attacker's expected payoff is at least as large at a solution of the RD method as at a CBSE of the ideal SG
\begin{equation}
    U^a(\bm{a}^*_o,\bm{d}^*_o) \leq  U^a(\bm{a}^*_{RD},\bm{d}^*_{RD}).
\end{equation}
\end{theorem}
\begin{proof}
 Please refer to Appendix C.
\end{proof}
To evaluate the defender's payoff loss due to its lack of information about the attacker's budget, we compute the mismatch  (loss) of actual defender's utility (\ref{eqactualpayoffrdsg}) of the RD method relative to that of CBSG (\ref{eq13b}) at its CBSE when the actual cost  is $\gamma_a$: 
\begin{equation}
    \mu_{\gamma_a} (\%) = \frac{U^d(\bm{a}^*_{RD},\bm{d}^*_{RD}) - U^d(\bm{a}^*_o,\bm{d}^*_o)}{U^d(\bm{a}^*_o,\bm{d}^*_o)} \times 100 \ \% \label{mismatchrdsg}
\end{equation}
Moreover, when the system realization is given by model $M_i \neq M_{nom}$, the proposed RD method developed for the nominal model can still be used, and the difference of the utilities obtained for the nominal model $M_{nom}$ and $M_i$ can be computed by (\ref{eq19}) with $(\bm{a}^*_{RD},\bm{d}^*_{RD})$ replacing $(\bm{a}^*_{nom},\bm{d}^*_{nom})$.

Finally, we note that the traversal method to compute the RD solution has at most $max(O((L_a+1)^n)$, $O((L_d+1)^n))$ computational complexity while the GA method to find this solution has $max(O((T+1)S_a)$,$O((T+1)S_d))$ complexity. Thus,  due to the sequential nature of RD algorithm, it is significantly simpler to implement than the ideal  CBSG (see Section III.D).

\section{Numerical Results for Wide-Area Control of Electric Power Systems}
\subsection{Power system model}
To demonstrate the performance of the proposed investment games, we consider one of the most important and safety-critical examples of an NCS, namely, an electric power system. The linear static state-feedback controller to be designed is referred to as wide-area control \cite{aranyapramod}, which helps in damping system-wide oscillations of power flows by minimizing the $\mathcal{H}_2$-performance function in (\ref{eq5}). Before discussing the game, we first briefly overview the dynamic model of the system. 
 
 Consider a power system network with $n$ synchronous generators and $\ell$ loads.
We assume that the state of the $i^{th}$ generator is denoted as $\bm{\xi}_i = [\delta_i, \; \omega_i, \; \bm{x}_{i,rem}] \in \mathbb R^{m_i}$ where $\delta_i$ is the generator phase, $\omega_i$ is the generator frequency, and $\bm{x}_{i,rem}$ is the vector of all non-electromechanical states.  The control input is considered as the excitation voltage and is denoted as $\Gamma_{i} \in \mathbb{R}$. All the loads in the system are considered as constant power loads without any dynamics. If one wishes to include dynamic loads, this can be easily accommodated by adding extra entries in the state vector $\bm{x}_i$. Let the pre-disturbance equilibrium of the $i^{th}$ generator be $\bm{\xi}^{0}_i$. The differential-algebraic model of the entire system, consisting of the generator models and the load models together with the power balance in the transmission lines,  is converted to a state-space model using Kron reduction (for details, please see \cite{Kundur1994}) and linearized about $\bm{\xi}_i^{0}, \,i = 1,2,...,n$.  The small-signal state of generator $i$ (or node $i$) is defined as $\bm{x}_i(t) = \bm{\xi}_i(t) - \bm{\xi}_i^{0}$. The small-signal model of the entire power system is thereafter written in the form of (\ref{eq1}). 

For the wide-area control design, $\bm{R}$ is chosen as the identity matrix while  $\bm{Q}=\mbox{diag}(\bar{\mathcal{L}},\,\bm{I}_n,\,\bm{I}_{m-2n})$ where $m = \sum_i^n m_i$, so that all generators arrive at a consensus in their small-signal changes in the phase angles. Here, $\bar{\mathcal{L}} = n\bm{I}_n - \bm{1}_n \cdot \bm{1}_n^T$  \cite{Lian2017JSAC} where $\bm{1}_n \in \mathbb{R}^{n \times 1}$ is the column vector of all ones, $\bm{I}_n$ is  the $n \times n$ identity matrix, and $\bm{I}_{m-2n}$ is  the $(m-2n) \times (m-2n)$ identity matrix.


\subsection{Case Study: IEEE 39-bus system}
 We employ the  IEEE 39-bus power model, which  consists of 10 synchronous generators and 19 loads, to evaluate the performance of the proposed SGs. Generator 1 is modeled by 7 states, generators 2 through 9 are modeled by 8 states each while generator 10 is modeled by 4 states. Therefore,  in this case $n = 9$. The dimension of the state vector in (\ref{eq1}) is $m = 75$, i.e. $\bm{A} \in \mathbb{R}^{75 \times 75}$, with the total number of control inputs $r = 9$, resulting in $\bm{B} \in \mathbb{R}^{75 \times 9}$ and $\bm{K} \in \mathbb{R}^{9 \times 75}$. We assume that $\bm{D} \in \mathbb{R}^{75 \times 9}$ in (\ref{eq1}) is a matrix with all elements zero except for the ones corresponding to the acceleration equation of all generators. The IEEE 39-bus system is referred to below as the {\it nominal model}. The data set for the IEEE 39-bus system along with the  randomly generated uncertain-model set $\mathcal{M}$ can be found in \cite{github_dataset}.

 \subsubsection{SG performance for fixed system model}
 \begin{table}[t]
\label{tabel1}
\large
\begin{center}
\resizebox{\columnwidth}{!}{%
\begin{tabular}[c]{ | P{2.3cm} | P{2.3cm} | P{2.5cm} | P{2.5cm} |} 
\hline
\rule{0pt}{10pt}
\bf{Node Rank} & \bf{Disabled generator} & \bf{Fract. loss $\%$  (local links disabled)} & \bf{Fract. loss $\%$ (local links intact}) \\[2pt]
\hline 
\hline 
1 & 9 & 40.7 &  0.07  \\[2pt]
\hline 
2 & 8 & 17.96 & 0.06 \\[2pt]
\hline 
3 & 4 & 17.67 & 0.06 \\[2pt]
\hline 
4 & 7 & 16.29 & 0.057\\[2pt]
\hline 
5 & 5 & 16.26 & 0.05 \\[2pt]
\hline 
6 & 3 & 14.63 & 0.045\\[2pt]
\hline
7 & 2 & 13.47 & 0.04\\[2pt]
\hline 
8 & 6 & 13.41 & 0.04 \\[2pt]
\hline 
9 & 1 & 3.79 & 0.01\\[2pt]
\hline
Open-loop & 1-9 & 181.92 & 0.4  \\[2pt]
\hline
\end{tabular}
}
\end{center}
\caption{Ranking of node ``importance" of the New England system (39-bus system) nominal model \cite{github_dataset} according to the fractional loss $\frac{\Delta_{\bm{s}_{m}}}{J(\bm{K}^{*}_{\mathcal{H}_2})} (\%)$ of the sparsity pattern $\bm{s}_{m}$ (\ref{eq9}) where only that node is attacked successfully. The open-loop fractional loss is also shown.}
\end{table}

  In subsections B.1 and B.2, the game is evaluated for the fixed nominal model available in \cite{github_dataset}.    Table \rom{1} illustrates the  fractional control performance losses $\frac{\Delta_{\bm{s}_{m}}}{J(\bm{K}^*_{\mathcal{H}_2})}(\%)$ (see (\ref{eq12})) for the  sparsity patterns where only one element of $\bm{s}_m$ in (\ref{eq9}) is zero, i.e. only one generator is  disabled. The highest loss is observed for the disabled generator $9$, followed by the  generators $8,4,7,5,3,2,6,1$, imposing the {\it``importance"} ranking of the generators, or nodes. From Table \rom{1}, we observe that when only the \textit{inter-node} communication links are disabled while the self-links are intact (see Fig. \ref{fig1}), the losses are greatly reduced, confirming that retaining the self-links is critical for a well-damped closed-loop performance. In other words, it is important to protect the system against the attacks that target the self-links \cite{PSACC2019}. The last row of Table \rom{1} shows the open-loop scenario with the loss $\Delta_{OL}$ (\ref{eq12l}). We observe very large loss in the latter case, especially when  the local links are disabled. However, note that the open-loop has finite $\mathcal{H}_2$-norm since in electric power systems, even if every link in the wide-area controller fails, the grid is still built to function stably by means of the power system stabilizers (PSS) inside the generators. For the rest of this section, we will make a practical assumption that a successful hardware or persistent DDoS attack on a given node disables both the self-links and the communication links to and from that node (see Fig. \ref{fig2}), reflecting the losses shown in the third column of Table \rom{1}.

\begin{figure}[htp]
    \centering
    \includegraphics[width=\linewidth]{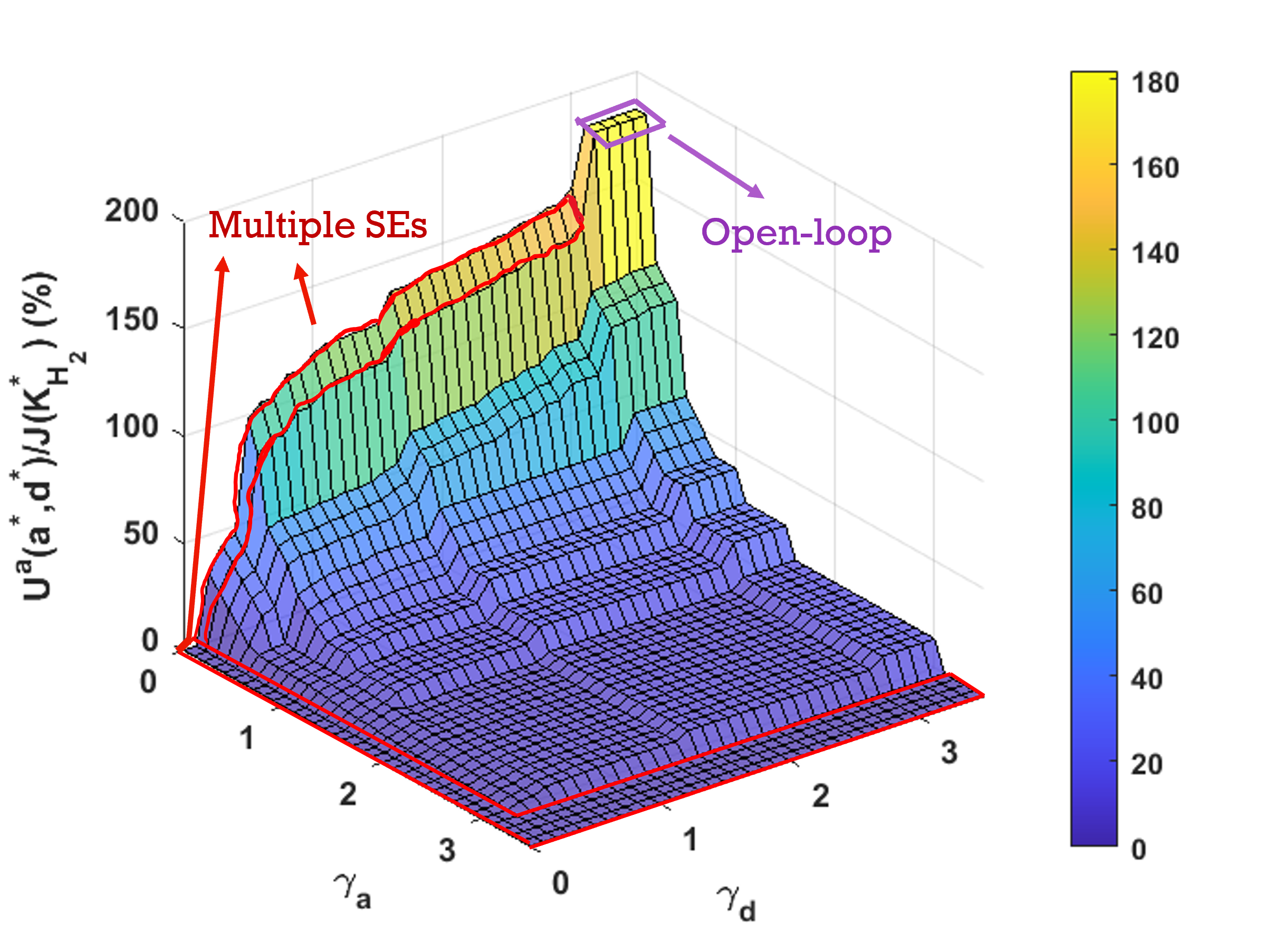}
    \caption{Attacker's fractional payoff at CBSE vs. costs per node $\gamma_a$ of the attacker  and $\gamma_d$ of the defender; $L_a = L_d = 3$; the 39-bus system nominal model \cite{github_dataset}. }
    \label{fig3}
\end{figure}

\begin{figure}[htp]
    \centering
    \includegraphics[width =\linewidth]{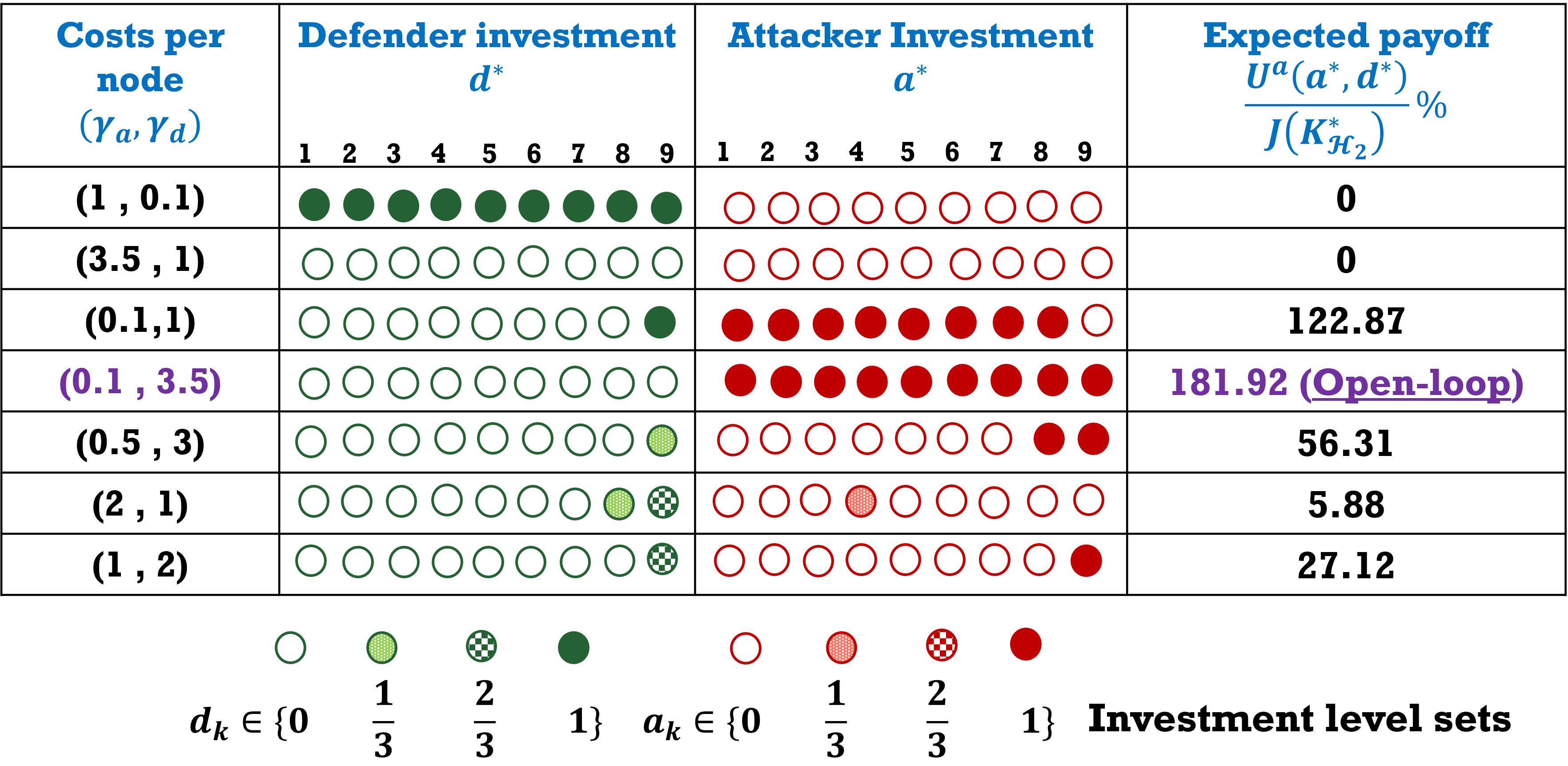}
    \caption{Players' strategies and the attacker's fractional payoff at CBSE for several ($\gamma_a,\gamma_d)$ pairs; $L_a=L_d=3$; the 39-bus system nominal model \cite{github_dataset}.}\label{fig4}
\end{figure}

Fig. \ref{fig3} shows the fractional attacker's payoff (\ref{eq13a}) (relative to $J(\bm{K}_{\mathcal{H}_2}^*)$) at CBSE versus the players' costs while Fig. \ref{fig4} illustrates the players' strategies and payoffs at CBSE. In both figures, we assume that all nodes of each player have the same cost per node in (\ref{eq11aa}), i.e., $\gamma_{a_i}=\gamma_a$, $\gamma_{d_i}=\gamma_d$ for all $i=1,...,n$. We observe that the performance trends are consistent with Theorem 1(a)$-$(d). 
In general, {\it multiple SEs} are possible for any choice of SG settings. 
In this example, they occur in the outlined regions of Fig. \ref{fig3}. The cost pair boundaries of these regions depend on the parameters $n=9$ and $L_a = L_d = 3$. First, multiple SEs  exist when $\gamma_d \leq \frac{1}{9}$, i.e., the defender is capable of protecting all nodes, resulting in $U^a(\bm{a}^*,\bm{d}^*) = 0$ (no loss) while 
the attacker chooses not to act at CBSE since it cannot change its payoff as shown in the first row of Fig. \ref{fig4}. The second region of multiple SEs corresponds to $\gamma_a > 3$ where the expected loss is zero since $\gamma_a > L_a$, i.e., the constraint (\ref{eq11aa}) is satisfied only when the attacker is inactive ($\bm{a} = \bm{0}$). Thus, the attacker does not have resources to attack in this region. In this case, a CBSE occurs when the defender also chooses not to  act (see the second row of Fig. \ref{fig4}).  Finally, multiple SEs exist when $\gamma_a \leq \frac{1}{9}$ while $\frac{1}{9} < \gamma_d \leq 3$ where the attacker is able to attack all nodes fully but both players choose cost-saving strategies. For example, in the third row of Fig. \ref{fig4}, the attacker saves its cost by not investing into node $9$ at CBSE since the defender fully protects this most ``important" node.

In the fourth row of Fig. \ref{fig4}, the defender is very {\it resource-limited} and thus does not act while the attacker attacks all nodes, resulting in the {\it open-loop} system, which occurs in the region $\gamma_a \leq \frac{1}{9}$, $\gamma_d > 3$ shown in Fig. \ref{fig3}. The last three rows of Fig. \ref{fig4} illustrate the scenarios where one or both players are resource-limited, and thus choose from the ``important" nodes (Table \rom{1}) to optimize their payoffs strategically. In the fifth row, the defender invests into the most ``important" node $9$ at the level $d_9 = \frac{1}{3}$ while the attacker has sufficient resources to invest fully into both ``important" nodes $9$ and $8$, thus raising the expected system loss above $50\%$. When $\gamma_a = 2$, $\gamma_d = 1$ (sixth row), the defender invests into  the ``important" nodes $9$ and $8$ while the attacker chooses the unprotected, but still ``important", node 4 since it has low chance of affecting the outcome for the more ``important" nodes due to its limited budget. The resulting attacker's payoff is low in this case.  Finally, in the last row, the defender is more resource-limited than the attacker, and both players target node 9 at the levels dictated by their cost constraints in (\ref{eq11aa}), with the resulting payoff increasing relative to the sixth row. 

While we assumed that all nodes have the same costs $\gamma_a$ or $\gamma_d$ for the attacker and defender, respectively, in some systems, the cost of attacking or protecting a certain node (e.g. an ``important" node) might be higher than that for other nodes. For example, we found that when a player's cost of the most ``important" node 9 increases significantly relative to the costs of the other nodes, that player avoids investing into node 9, and thus it's payoff decreases.

\begin{figure}[ht]
  \begin{center}
    \begin{minipage}{0.48\textwidth}
      \includegraphics[width=\linewidth]{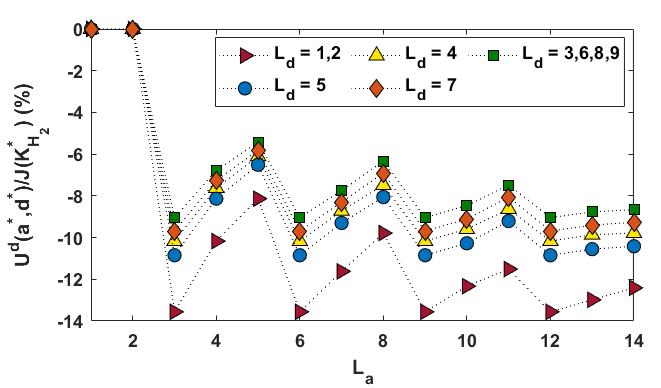}
      \caption*{(a)}
       \end{minipage}
       \begin{minipage}{0.48\textwidth}
       \includegraphics[width=\linewidth]{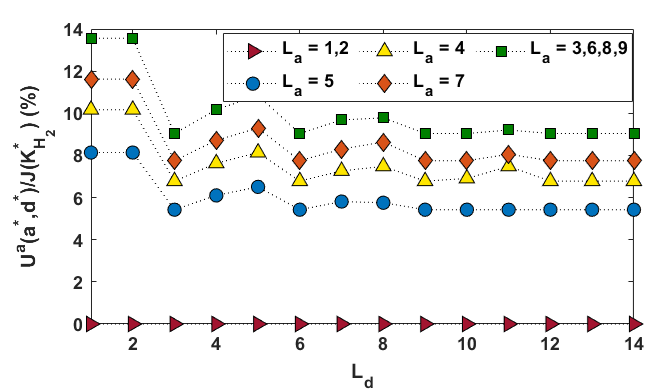}
       \caption*{(b)} 
      \end{minipage}
\caption{Players' fractional payoffs at CBSE;  (a) Defender's fractional payoff  vs. $L_a$ for several $L_d$ values; (b)  Attacker's fractional payoff  vs. $L_d$ for several $L_a$ values; $\gamma_a=\gamma_d=1.5$; the 39-bus system nominal model \cite{github_dataset}.}\label{fig5}
   \end{center}
\end{figure}
In Fig. \ref{fig5}, we illustrate the players' fractional payoffs for $\gamma_a=\gamma_d = 1.5$ (which characterize moderate resources of both players) as a player's number of investment levels varies while its opponent's number of investment levels is fixed. Similar simulations were performed for other cost pairs, and the results are consistent with Theorem 1(e). 
We observe that setting a player's number of levels to {\it three} provides a near-optimal payoff for that player for any fixed opponent's number of levels. Thus, $L_a = L_d = 3$ is used in the numerical results throughout the paper.  

Finally, in \cite{PS_thesis} we compared the proposed game with {\it individual optimization} (IO) where the attacker aims to degrade the system performance by increasing the $\mathcal{H}_2$-performance cost $J$ in (\ref{eq5}) while the defender aims to decrease it. We summarize these results below. In IO, both players act  under their individual budget constraint and without taking into account the opponent's possible actions. Thus, as in the CBSG, each player always targets its ``important" nodes in IO.  However, in the game the players' strategies affect each other and the attacker often prefers to avoid the nodes protected by the defender as  was illustrated in Fig. \ref{fig4}. This behavior was not observed in IO. We found that each player's payoff can be up to $50\%$ lower when using IO relative to playing the game for some cost pairs $(\gamma_a,\gamma_d)$. Moreover, the players' costs can be  significantly lower when playing the proposed CBSG than in IO since the CBBI Algorithm \ref{alg:ALG1} selects an SE with the lowest cost while in IO each player invests fully up to its budget constraint \cite{PS_thesis}.

\begin{figure}[ht]
  \begin{center}
    \begin{minipage}{0.48\textwidth}
      \includegraphics[width=\linewidth]{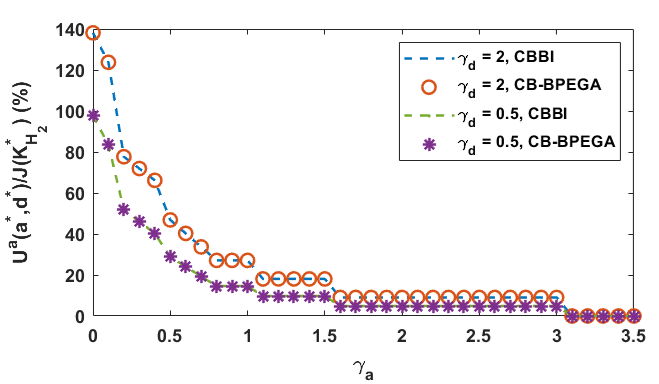}
      \caption*{(a)}
      \end{minipage}
            \begin{minipage}{0.48\textwidth}
      \includegraphics[width=\linewidth]{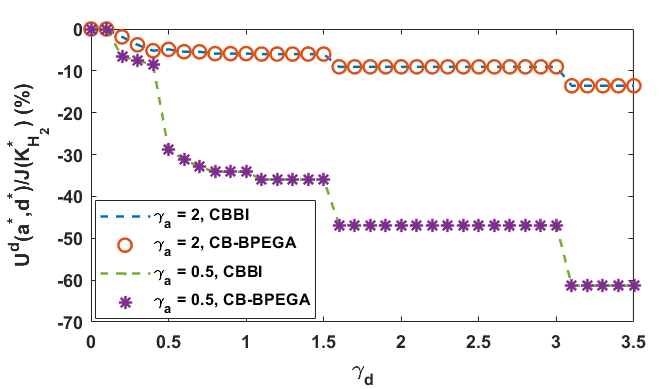}
      \caption*{(b)}
      \end{minipage}
\caption{ Players’ payoffs at the CBSE found by the CBBI Algorithm \ref{alg:ALG1} and the CB-BPEGA Algorithm \ref{alg:ALG2}
at convergence: (a) Attacker’s payoff vs. $\gamma_a$; 
(b) Defender’s payoff vs. $\gamma_d$; $L_a = L_d = 3$; the 39-bus system nominal model~\cite{github_dataset}.}\label{GAfig7}
   \end{center}
\end{figure}

\subsubsection{CB-BPEGA validation}
In this subsection, we validate convergence of the CB-BPEGA Algorithm \ref{alg:ALG2} for the  nominal model of the IEEE 39-bus system. We set the attacker's population size $S_a = 30$ and the defender's population size $S_d = 30$. The crossover probability $P_c = 0.85$ and the mutation rate $P_m = 0.10$. The maximum number of generations is set to $T = 20$. These initialization parameters are selected experimentally and reflect the trade-off between the convergence and computational complexity for the IEEE 39-bus system. The CB-BPEGA method described in Algorithm \ref{alg:ALG2} is applied to find a CBSE for the proposed CBSG.

\begin{figure}[htp]
    \centering
    \includegraphics[width=\linewidth]{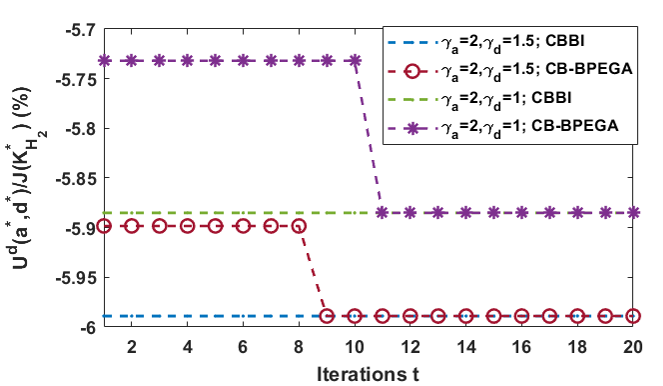}
    \caption{Convergence of the defender’s payoff in the CB-BPEGA Algorithm \ref{alg:ALG2} with $L_a=L_d=3$; the 39-bus system nominal model \cite{github_dataset}.}\label{GAfig8}
\end{figure}

Fig. \ref{GAfig7} shows the comparison of the players' payoffs at a CBSE of the CBBI Algorithm \ref{alg:ALG1} and of the CB-BPEGA Algorithm \ref{alg:ALG2} at convergence versus cost of attack $\gamma_a$ (Fig. \ref{GAfig7}a) and cost of defense $\gamma_d$ (Fig. \ref{GAfig7}b) for fixed levels of investment $L_a = L_d = 3$. We observe that the CB-BPEGA method converges to  a CBSE found by the CBBI Algorithm \ref{alg:ALG1}, confirming Proposition 5.1 in \cite{Lu_thesis}.
 Moreover, Fig. \ref{GAfig8} demonstrates that Algorithm \ref{alg:ALG2} converges to a CBSE obtained by the traversal CBBI method in  fewer than 15 iterations. Similar results were obtained for other cost pairs, demonstrating fast convergence of the CB-BPEGA Algorithm \ref{alg:ALG2}.

For the nominal model \cite{github_dataset} of the IEEE 39-bus system, given  parameter values above,  CBBI Algorithm \ref{alg:ALG1} has at most as $O(4^9 \times 4^9)$ complexity from Remark 3 while the computational complexity of the CB-BPEGA Algorithm \ref{alg:ALG2} is $O(20 \times 30 \times 30) = O(18000)$ for any cost pair. 

\subsubsection{Performance of the nominal-model CBSG for uncertain systems}
In this subsection,  to evaluate the robustness of the proposed CBSG (Section IV.A) to model uncertainties, we generate a set of uncertain models of the IEEE 39-bus system. This is done by perturbing the reactive power setpoints of all loads and the  inertias of all generators by adding independent and identically distributed (iid) Gaussian random variables  \cite{billinton2008effects}
with {\it zero} means ($\mu = 0$) and {\it unity} standard deviations ($\sigma = 1$). 
The probability of the change in the load reactive power is considered to be twice that of the change in inertia \cite{Kundur1994}. For each combination of the perturbations in these parameters, power flow is run on the nonlinear power system model to determine the corresponding equilibrium, followed by small-signal linearization around this equilibrium. The model set is represented as $\mathcal{M} = \left\{M_1,...,M_j,...,M_M\right\}$ where $M_j = \left\{\bm{A}_j,\bm{B}\right\}$. 
As the matrices $\bm{B}$ and $\bm{D}$ in (\ref{eq1}) are independent of the chosen uncertainties,  we assume them to be fixed at the nominal values for all models in $\mathcal{M}$. 
Note that the number of models ($M$) in the set $\mathcal{M}$ satisfies the criteria for the sample size given the margin of error (MOE) of $0.05$, the standard deviation of $1$, and the confidence interval (CI) of $95\%$. The sample size or the number of  models in the set $\mathcal{M}$ is computed as $M \geq  ((z \times \sigma)/ MOE)^2 = 1536.64$ where $z = 1.96$. Thus, we set $M = 1550$.
This sample size takes into account uncertainties in both loads and generator inertia~\cite{HowToStat,kotz1999encyclopedia}.

\begin{figure}[htp]
    \centering
    \includegraphics[width=\linewidth]{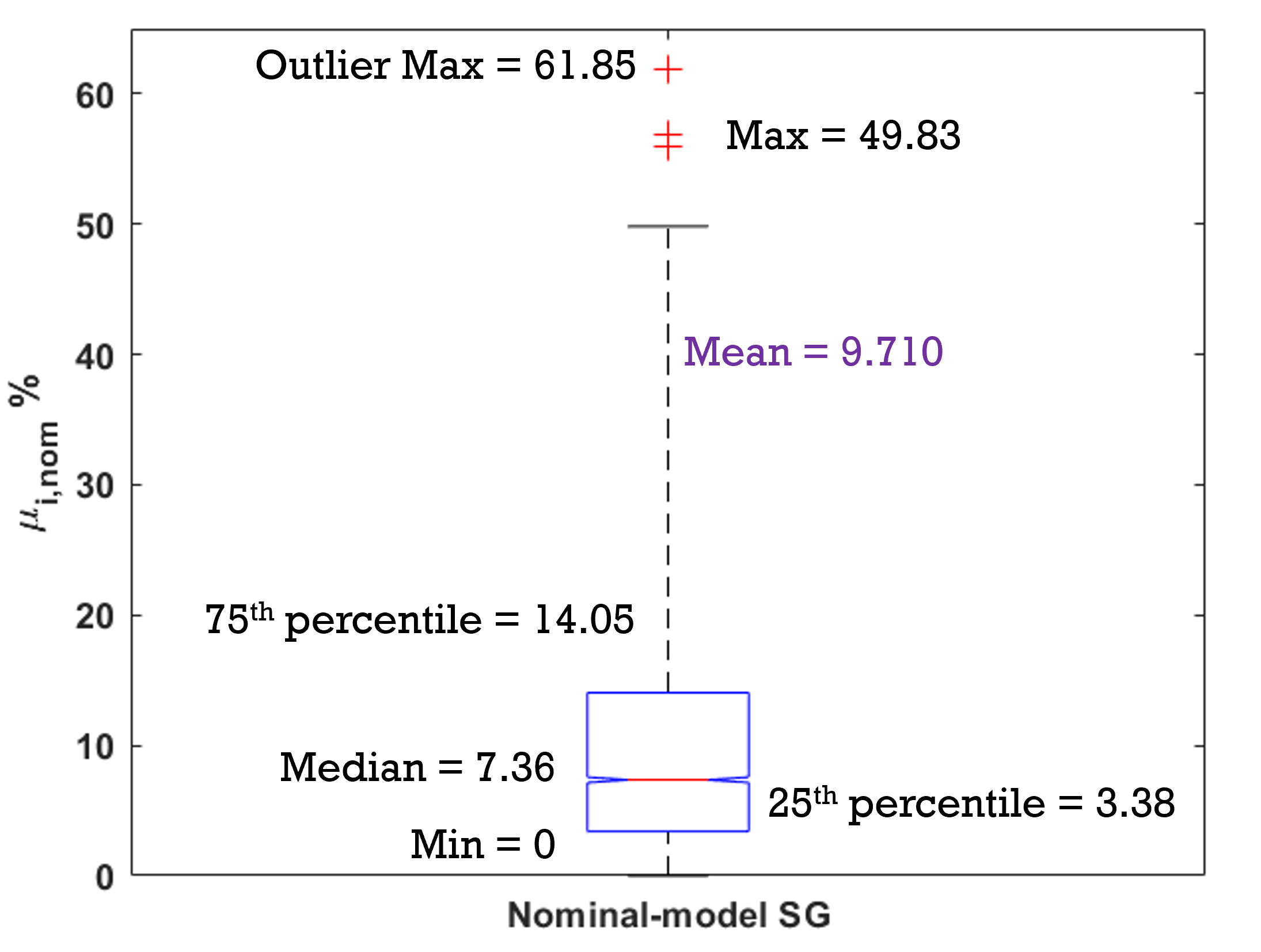}
    \caption{Boxplot of the utility difference $\mu_{i,nom}\%$ (\ref{eq19})  over a randomly generated set $\mathcal{M}$ of  $1550$ models  and $4$ cost pairs $(\gamma_a,\gamma_d)$: $(1.5,1.5)$, $(0.5,1.5)$, $(1.5,0.5)$ and  $(2.5,2.5)$; $L_a = L_d = 3$; the 39-bus system~\cite{github_dataset}.}\label{fig12}
\end{figure}

Fig. \ref{fig12} represents the {\it utility difference} (\ref{eq19}) statistics for the nominal-model SG employed over  the set $\mathcal{M}$ for certain cost pairs. Since the nominal-payoff model $M_1$ is included in $\mathcal{M}$, $\mu_{1,nom}=\mathop { \min }\limits_{i} \mu_{i,nom}=0$. We found that for individual cost pairs within Fig. \ref{fig12}, these statistics were similar. 
For example, the mean mismatch values for cost pairs $(1.5,1.5)$, $(0.5,1.5)$, $(1.5,0.5)$, and $(2.5,2.5)$ are $9.55$, $10.17$, $11.67$ and $12.78$, i.e. within close proximity of each other \cite{PS_thesis}. We note that most uncertain models experience modest payoff differences from the nominal-model payoff shown in Fig. \ref{fig3}. Thus, the nominal-model CBSG is a robust solution to security investment in~NCS.

\subsubsection{ Performance of Robust-Defense Method} 
\begin{figure}[htp]
  \begin{center}
    \begin{minipage}{0.49\textwidth}
      \includegraphics[width=\linewidth]{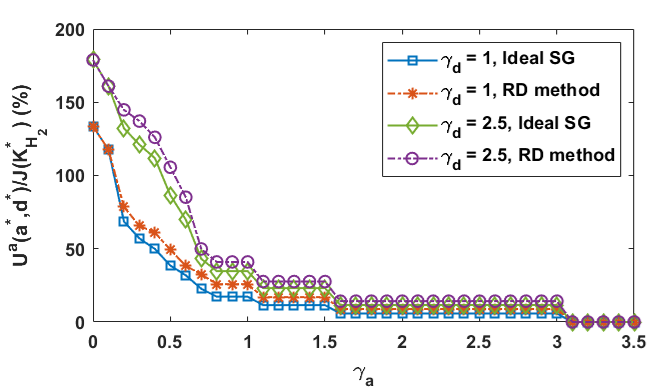}
      \caption*{(a)}
       \end{minipage}
       \begin{minipage}{0.49\textwidth}
       \includegraphics[width=\linewidth]{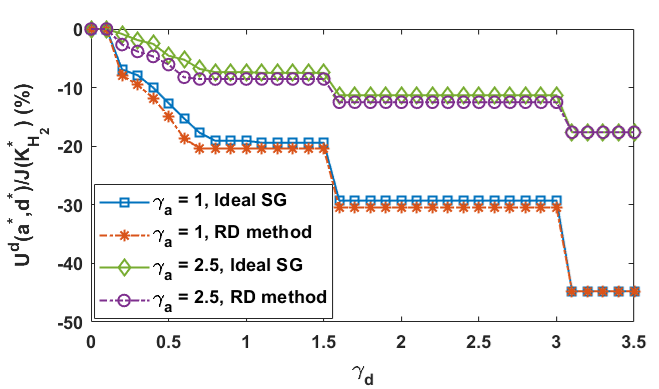}
       \caption*{(b)} 
      \end{minipage}
\caption{Players' actual fractional payoffs at CBSE of the nominal-model SG (Section III.C) and the RD solution (Section IV);  (a) Attacker's fractional expected payoff  vs. $\gamma_a$ for $\gamma_d = 1,2.5$; (b) Defender's actual fractional payoff  vs. $\gamma_d$ for $\gamma_a = 1,2.5$; $L_a = L_d = 3$; the 39-bus system \cite{github_dataset}.}\label{fig13}
   \end{center}
\end{figure}

In this subsection, we compare the performance of the RD Algorithm \ref{alg:ALG3} and the ideal SG Algorithm \ref{alg:ALG1} when the system model is fixed as the nominal model.  Fig. \ref{fig13} shows the fractional payoffs of the players for both methods for selected costs of attack and defense. These results are consistent with Theorem 2. Moreover, we note that the players' payoffs for the realistic RD method are close to those of the ideal SG, which assumes complete knowledge of the attacker's resources by the defender.

\begin{figure}[ht]
    \centering
    \includegraphics[width=\linewidth]{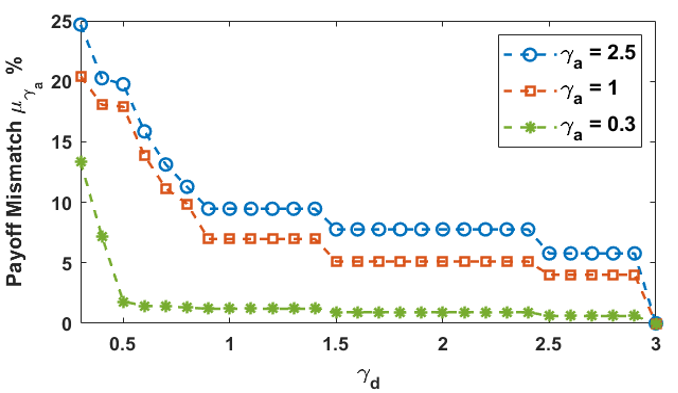}
    \caption{Mismatch $\mu_{\gamma_a}$ ($\%$) (\ref{mismatchrdsg}) of the RD method due to uncertain attacker's resources for the defender when the actual attacker's cost is $\gamma_a$; the 39-bus system \cite{github_dataset}.}\label{fig13new}
\end{figure}

For the cost pair values in Fig. \ref{fig13}, Fig. \ref{fig13new} shows the {\it mismatch} (loss) of the defender's payoff in the RD method  $\mu_{\gamma_a}$ ($\%$) (\ref{mismatchrdsg}) due to the defender's lack of knowledge of the attacker's budget. 
We observe that the mismatch can be as high as $24.9\%$ (for $\gamma_a = 2.5, \gamma_d=0.3$) due to overprotection, but  tends to decrease as $\gamma_d$ grows and $\gamma_a$ decreases. The former trend is due to limited protection options for a resource-constrained defender while the latter is caused by a better match between the actual $\gamma_a$ value and the defender's estimate of zero attacker's cost.

\begin{figure}[ht]
    \centering
    \includegraphics[width =\linewidth]{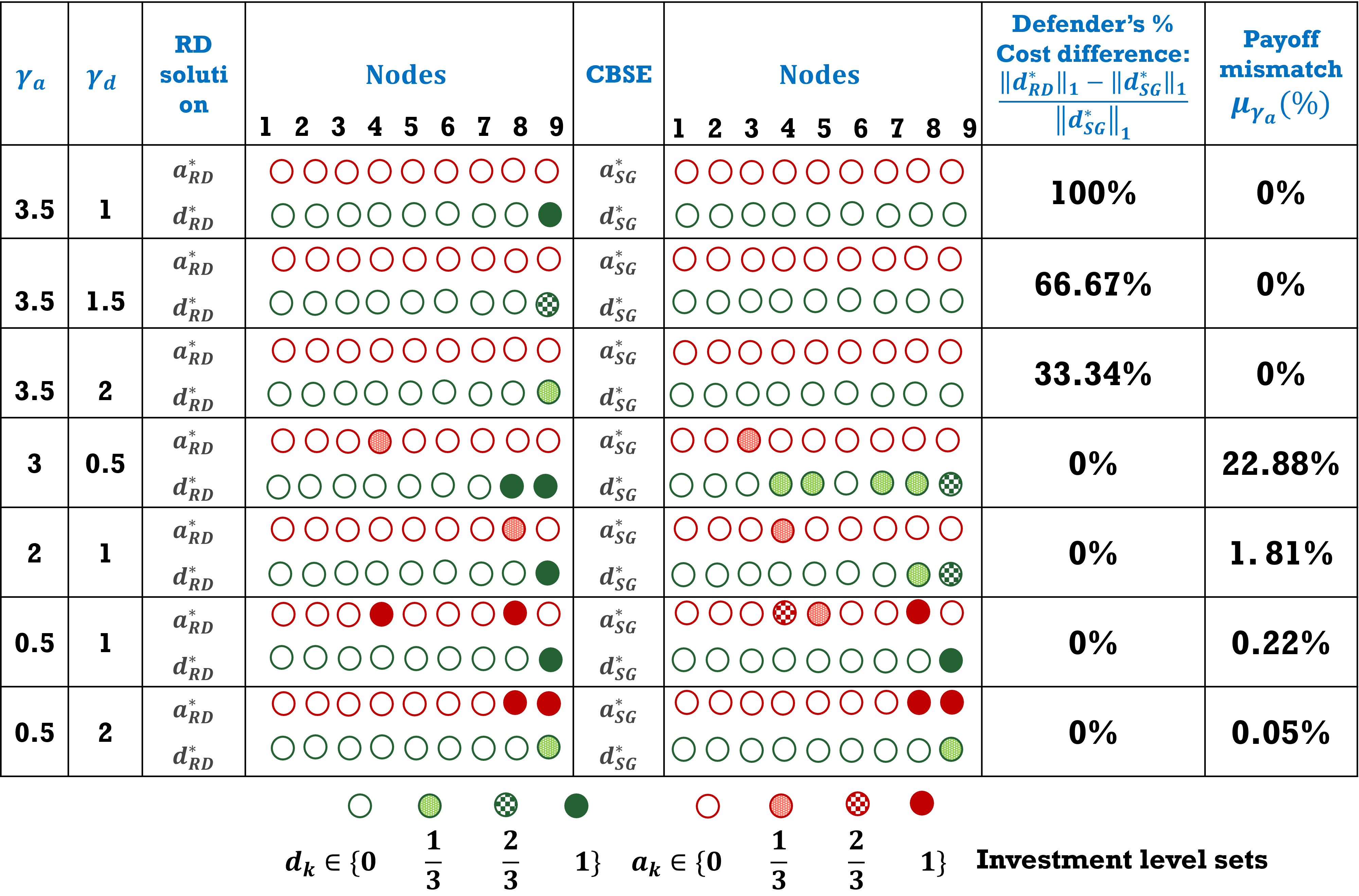}
    \caption{Players' strategies and defender's excess cost and mismatch (loss) (\ref{mismatchrdsg}) at  the solution of the RD method $(\bm{a}^*_{RD},\bm{d}^*_{RD})$ (\ref{eq:ch4.32})-(\ref{eq:ch4.35}) vs. ideal SG at CBSE $(\bm{a}^*_{SG},\bm{d}^*_{SG})$ (\ref{eq18}) for several ($\gamma_a,\gamma_d)$ pairs; $L_a=L_d=3$; the 39-bus system \cite{github_dataset}.}\label{fig14}
\end{figure}

 Furthermore, due to the assumption of the most powerful attacker, the defender always invests fully up to the given budget constraints (\ref{eq11aa}) and thus is likely to overinvest (incur a higher cost) in the RD method relative to the ideal SG where  the defender knows the attacker's resources. On the other hand, the attacker's investment cost is  the same in both methods. Fig. \ref{fig14} shows the  players' strategies and the defender's excess cost of the RD method relative to the ideal CBSG for several pairs of the actual costs of attack/defense per node. The figure also shows the  payoff mismatch (\ref{mismatchrdsg}) for these cost pairs.  
Note that when the cost of attack per node is high (e.g., $\gamma_a = 3.5$), the defender's cost is much higher in the RD  method than in the ideal SG since  in the former the defender overestimates the attacker's resources and thus commits more resources than necessary for protection. We note that the actual payoffs of both methods are zero in this case due to the weak attacker. When the attacker is stronger ($\gamma_a \leq 3$), the defender's costs are the same in both approaches but the defender is able to commit its resources more strategically in the ideal SG since it knows the attacker's resources. For example, in the  {\it fourth} row of Fig. \ref{fig14}, the defender invests fully in the two most important nodes in the RD method but distributes its resources over $5$ nodes in the  ideal SG. The resulting utility of the defender is $22.88\%$ lower in the RD method. In the  {\it fifth} row of Fig. \ref{fig14}, the defender's utility loss relative to the ideal SG is much smaller ($1.81\%$) since the attacker is relatively strong in this case. Similarly, in the last two rows, the mismatch further reduces due the attacker's low cost ($\gamma_a = 0.5$), thus approximating the assumption of its zero cost made by the defender in the RD algorithm. 

Finally, we found that the RD GA method in Section IV.B is able to compute the RD solution in about $15$ iterations. The parameters of this method were the same as that for the CB-BPEGA method in Section V.B.2. We also simulated the RD method for the uncertain model set $\mathcal{M}$ and observed that the payoff difference (\ref{eq19}) and boxplot statistics were similar to the results for the nominal-model SG in Section V.B.3 \cite{PS_thesis}.

\subsubsection{Computational Complexity}
 This subsection summarizes the computational requirements of the proposed methods for the IEEE 39-bus model. First, all algorithms require the computation of the loss vector corresponding to all sparsity patterns in (\ref{eq9}) (see Remark 3). This computation is  dominated by the structured $\mathcal{H}_2$ optimization algorithm in \cite{Mihailo2013TAC}, which is completed in under $900$ seconds\footnote{The experiments are run using MATLAB on Windows 10 with 64-bit operating system, 3.4 GHz Intel core i7 processor, and 16GB memory. \label{note1}}.
 
 For Algorithm \ref{alg:ALG1} (CBBI method), the payoff matrix (\ref{eq13a}) must be computed for all actions that satisfy (\ref{eq11aa}). This computation scales with the number of feasible actions or payoffs (Remark 3), which is bounded by $4^{18}$ for the SG designed for the IEEE 39-bus system model with $L_a=L_d =3$. On the other hand, for Algorithm \ref{alg:ALG2} (CB-BPEGA), the entire matrix does not need to be computed. Instead, the algorithm scales with the population sizes $S_a$, $S_d$, reducing the size of the payoff matrix.
We found that CB-BPEGA Algorithm \ref{alg:ALG2} reduced the time complexity significantly (see Section III.D) for lower cost values, i.e. larger action space.  For example, when $\gamma_a = \gamma_d = 0.5$, the action space of the players is not restricted significantly,  corresponding to  a high computation load. In this case, the computation time of CBSE is about $1120$ seconds using CBBI Algorithm \ref{alg:ALG1}  and is under $410$ seconds for the CB-BPEGA Algorithm \ref{alg:ALG2} with the number of iterations $T=20$ and population sizes as $S_a=S_d=30$. 

Moreover,  as discussed in Section IV.B, the RD method has lower computation load than the CBSG.
When $\gamma_a=\gamma_d=0.5$, the  computation time of the  Algorithm \ref{alg:ALG3} (sequential traversal method)  is $460$ seconds while the GA for finding the RD solution with the number of iterations $T=20$ and population sizes as $S_a=S_d=30$ is completed in under $150$ seconds. Thus, the latter reduces the computational complexity by a factor of $2.5$  relative to the former.

Finally, we note that the proposed methods solve a long-term resource-planning problem and are computed offline, so the computation complexity does not significantly impact their implementation.

\subsection{Extension to IEEE-68 bus power system model}
The 68-bus test system model \cite{Book-RobustControl} consists of 16 generators and 52 load buses. The generators are modeled by 10 states each. Thus, in (\ref{eq1}), the number of nodes in the NCS is $n=16$, the dimension of the state vector is $m=160$, i.e. $\bm{A}\in \mathbb{R}^{160 \times 160}$, the total number of control inputs $r=16$, resulting in $\bm{B}\in \mathbb{R}^{160 \times 16}$ and $\bm{K} \in \mathbb{R}^{16 \times 160}$. We assume that $\bm{D} \in \mathbb{R}^{160 \times 16}$ is a matrix with all elements equal to zero except for those corresponding to the acceleration equation of all generators. 

\begin{table}[ht]
    \centering
    \large
    \resizebox{\columnwidth}{!}{%
    \begin{tabular}{|P{1 cm} | P {1.5 cm} | P {1 cm} || P {1cm} | P {1.5cm} | P {1cm} |} 
    \hline
         \bf{Node rank} & \bf{Disabled gen} & \bf{Fract. loss $\%$}  & \bf{Node rank} & \bf{Disabled gen} & \bf{Fract. loss $\%$ } \\
         \hline \hline
         1 & 6 & 19.66 & 9 & 14 & 19.08\\
         \hline
          2 & 8 & 19.61 & 10 & 1 & 19.06 \\
        \hline 
        3 & 13 & 19.32 & 11 & 2 & 19.02 \\
        \hline 
        4 & 4 & 19.30 & 12 & 15 & 19.01 \\
        \hline 
        5 & 9 & 19.20 & 13 & 11 & 18.86 \\
        \hline 
        8 & 10 & 19.19 & 14 & 5 & 18.83\\
        \hline
        7 & 12 & 19.15 & 15 & 7 & 18.72\\
        \hline 
        8 & 3 & 19.14 & 16 & 16 & 8.94\\
        \hline 
        \multicolumn{6}{|c|}{Open-loop fract. loss ($\%$) - 169.87}\\
        \hline
    \end{tabular}
    }
    \caption{Ranking of node ``importance" according to the fractional loss $\frac{\Delta_{\bm{s}_{m}}}{J(\bm{K}^{*}_{\mathcal{H}_2})} (\%)$ of the sparsity pattern $\bm{s}_{m}$ (\ref{eq9}) where only that node is attacked successfully. The open-loop fractional loss is $169.87$ for disabled generators $1-16$; the IEEE 68-bus system nominal model \cite{github_dataset}. }
    \label{table3}
\end{table}
Table \ref{table3} illustrates the fractional control performance losses $\frac{\Delta_{\bm{s}_{m}}}{J(\bm{K}^*_{\mathcal{H}_2})}(\%)$ (see (\ref{eq12})) for the  sparsity patterns where only one element of $\bm{s}_m$ in (\ref{eq9}) is zero, i.e. only one generator is  disabled for the ``nominal" model of the IEEE 68-bus system.. The highest loss is observed for the disabled generator $6$, followed by the  generators $8,13,4$ and so on, imposing the {\it``importance"} ranking of the generators. Note that this fractional loss represents a successful hardware or persistent DDoS attack on a node that disables both self-links and the communications links to and from that node (see Fig. \ref{fig2}).

\begin{figure}[ht]
  \begin{center}
    \begin{minipage}{0.49\textwidth}
      \includegraphics[width=\linewidth]{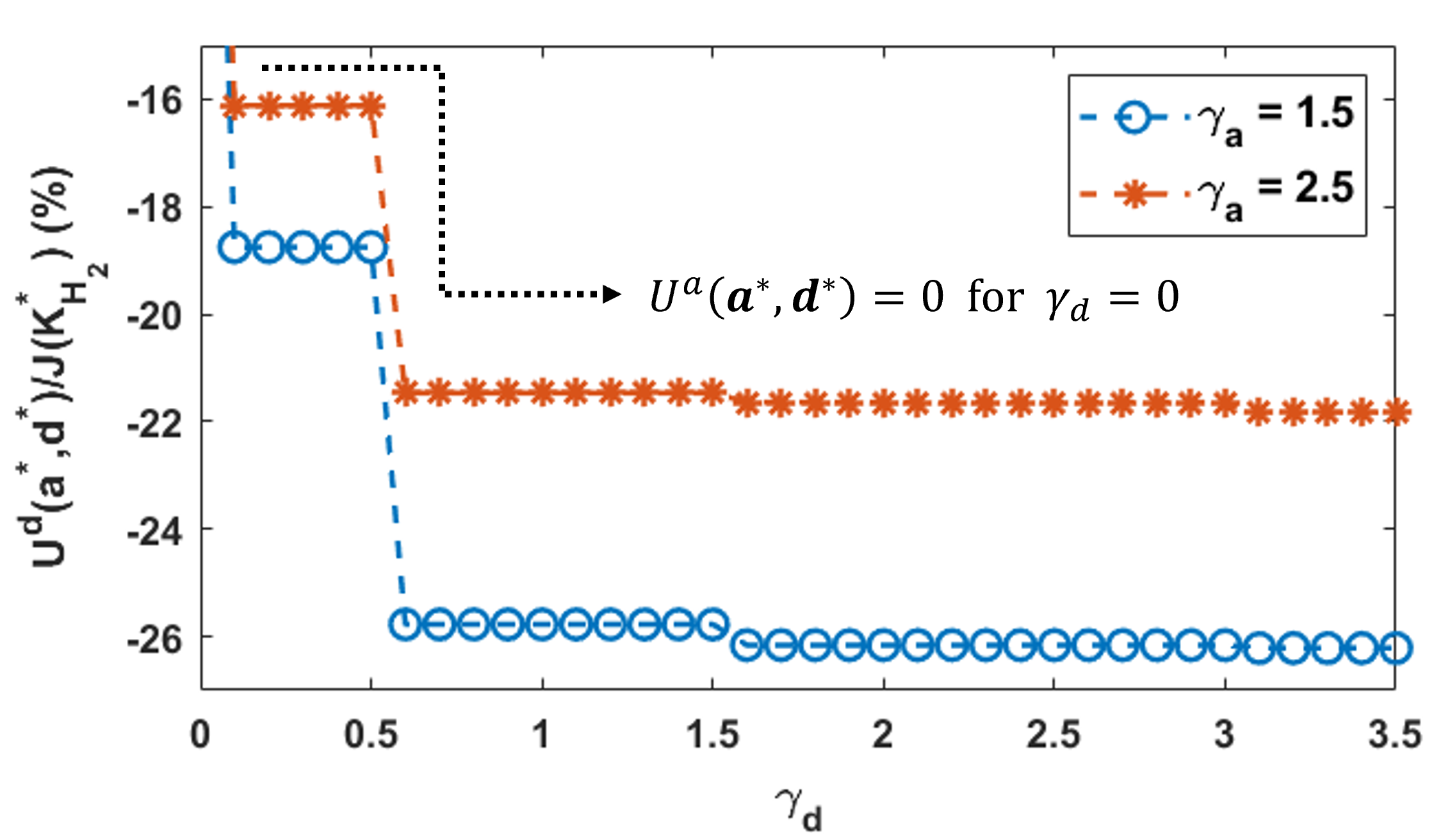}
      \caption*{(a)}
       \end{minipage}
       \begin{minipage}{0.46\textwidth}
       \includegraphics[width=\linewidth]{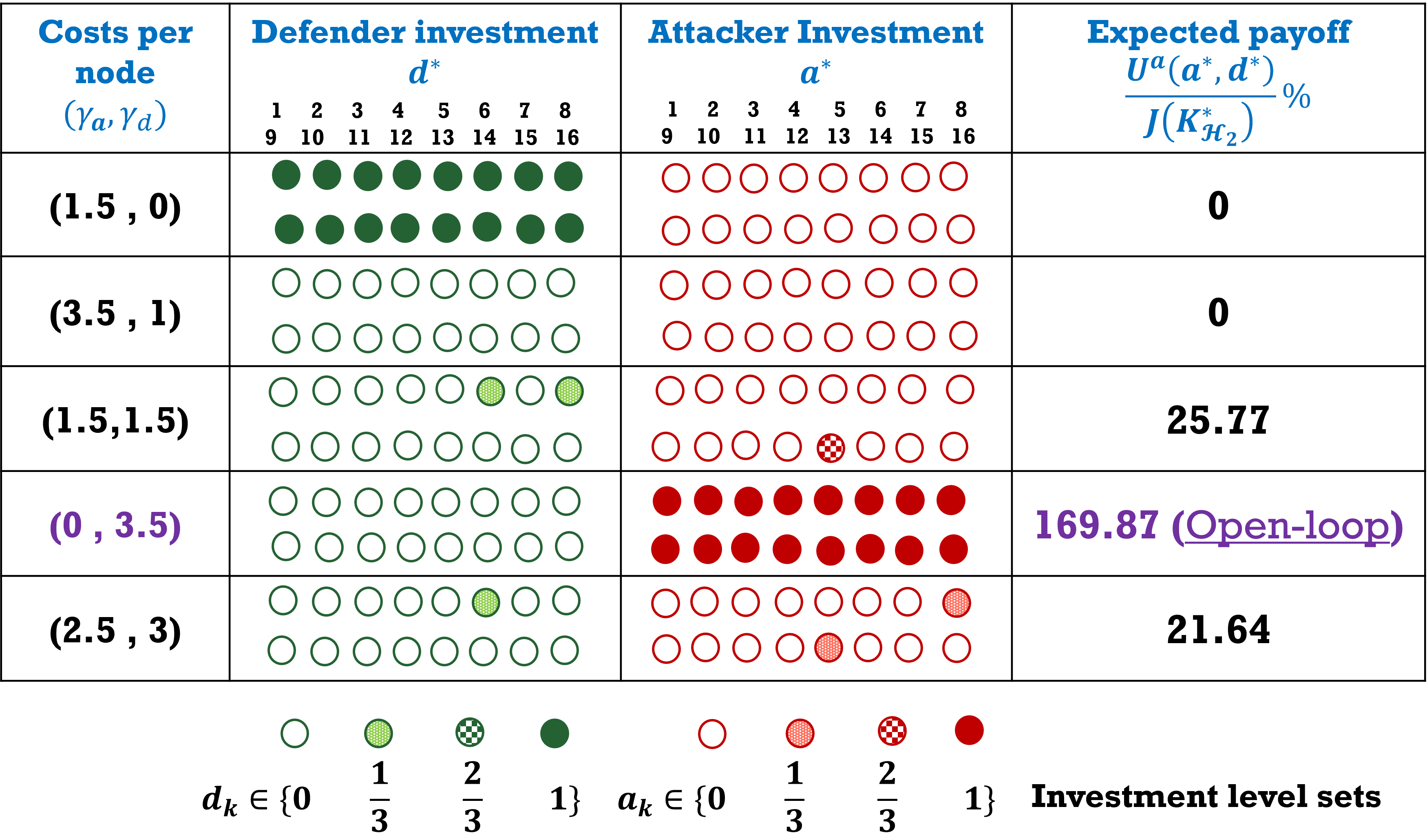}
       \caption*{(b)} 
      \end{minipage}
\caption{(a) Defender's fractional payoffs at CBSEs vs. $\gamma_d$ for $\gamma_a = 1.5,2.5$;  (b) Players' strategies and the attacker's fractional payoff at CBSE for several ($\gamma_a,\gamma_d)$ pairs; $L_a=L_d=3$; the 68-bus system nominal model \cite{github_dataset}.}\label{fig15}
   \end{center}
\end{figure}
The CBBI Algorithm \ref{alg:ALG1} and traversal sequential optimization Algorithm \ref{alg:ALG3}  have very high complexity  for the 68-bus model, so  in this subsection we employ only the CBPEGA Algorithm \ref{alg:ALG2} and the sequential GA method (Section IV.B). For both GA methods, we set the attacker's population size $S_a = 20$ and the defender's population size $S_d = 30$. The crossover probability $P_c = 0.85$, and the mutation rate $P_m = 0.10$. The maximum number of generations is set to $T = 25$. These initialization parameters are selected experimentally and reflect the trade-off between the convergence and computational complexity for the IEEE 68-bus system. 

Assuming the fixed, nominal model, Fig. \ref{fig15}a shows the fractional defender's payoff (\ref{eq13b}) (relative to $J(\bm{K}_{\mathcal{H}_2}^*)$) at CBSE versus the cost of defense for two values of attacker's costs while Fig. \ref{fig15}b illustrates the players' strategies and payoffs at CBSE  for several $(\gamma_a,\gamma_d)$ cost pairs. In both figures, we assume that all nodes of each player have the same cost per node in (\ref{eq11aa}). We observe that the performance trends are consistent with Theorem 1(a)$-$(d). 

Next,  we evaluate the performance of the nominal-payoff CBSG over uncertain models in set $\mathcal{M}$ of the IEEE 68-bus system.  Similarly to the IEEE 39-bus system,  we generate a set of uncertain models for the IEEE 68-bus system by perturbing the reactive power setpoints of every load and the inertia of every generator by adding iid Gaussain variables with zero mean and unity variance to their nominal values. Note that the number of models ($M$) in the set $\mathcal{M}$ is chosen as $1550$ to satisfy the statistical sample size as in Section V.B.3. 
\begin{figure}[ht]
    \centering
    \includegraphics[width=\linewidth]{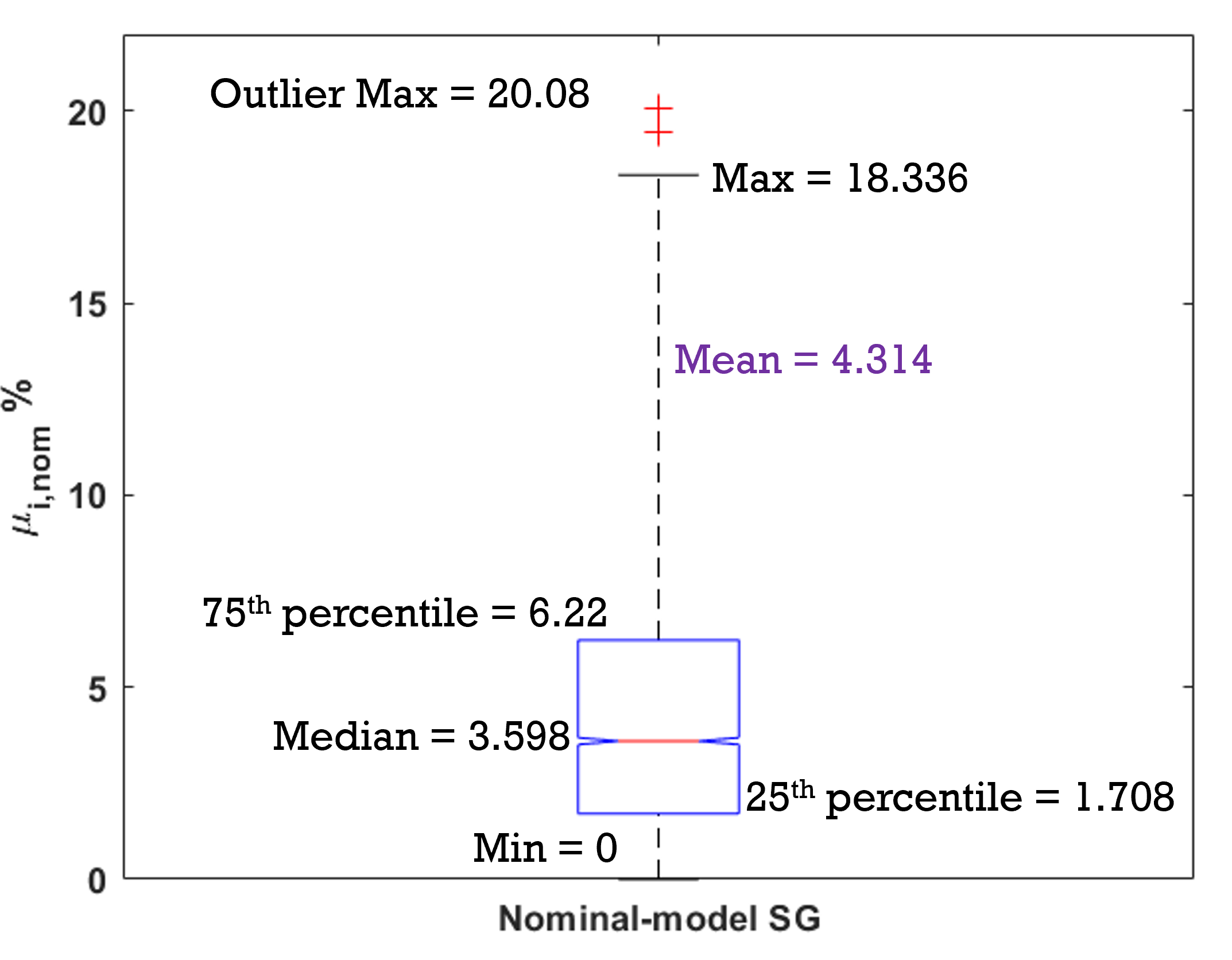}
    \caption{Boxplot of the utility difference $\mu_{i,nom}\%$ (\ref{eq19})  over a randomly generated set $\mathcal{M}$ of  $1550$ models  and $4$ cost pairs $(\gamma_a,\gamma_d)$: $(1.5,1.5)$, $(2.5,2.5)$, $(1.5,2.5)$ and  $(2.5,1.5)$; $L_a = L_d = 3$; the 68-bus system \cite{github_dataset}.}\label{fig16}
\end{figure}
Fig. \ref{fig16} shows the {\it payoff difference} (\ref{eq19}) statistics for the nominal-model SG employed over $1550$ models in $\mathcal{M}$ for certain cost pairs. Similarly to the results for the 39-bus model in Fig. \ref{fig12}, 
we observe that the nominal-model game is a robust approach to security investment in  this large-scale system.

\begin{figure}[h]
  \begin{center}
    \begin{minipage}{0.48\textwidth}
      \includegraphics[width=\linewidth]{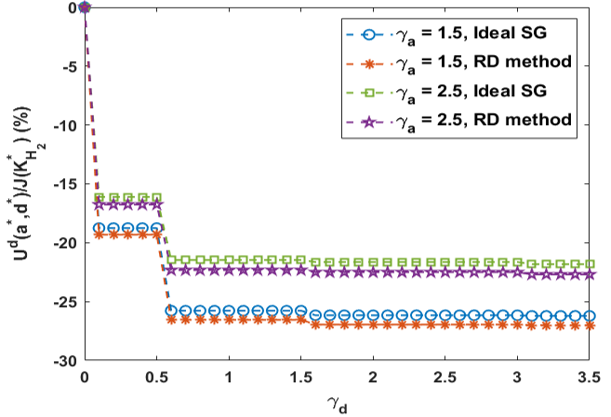}
      \caption*{(a)}
       \end{minipage}
       \begin{minipage}{0.465\textwidth}
       \includegraphics[width=\linewidth]{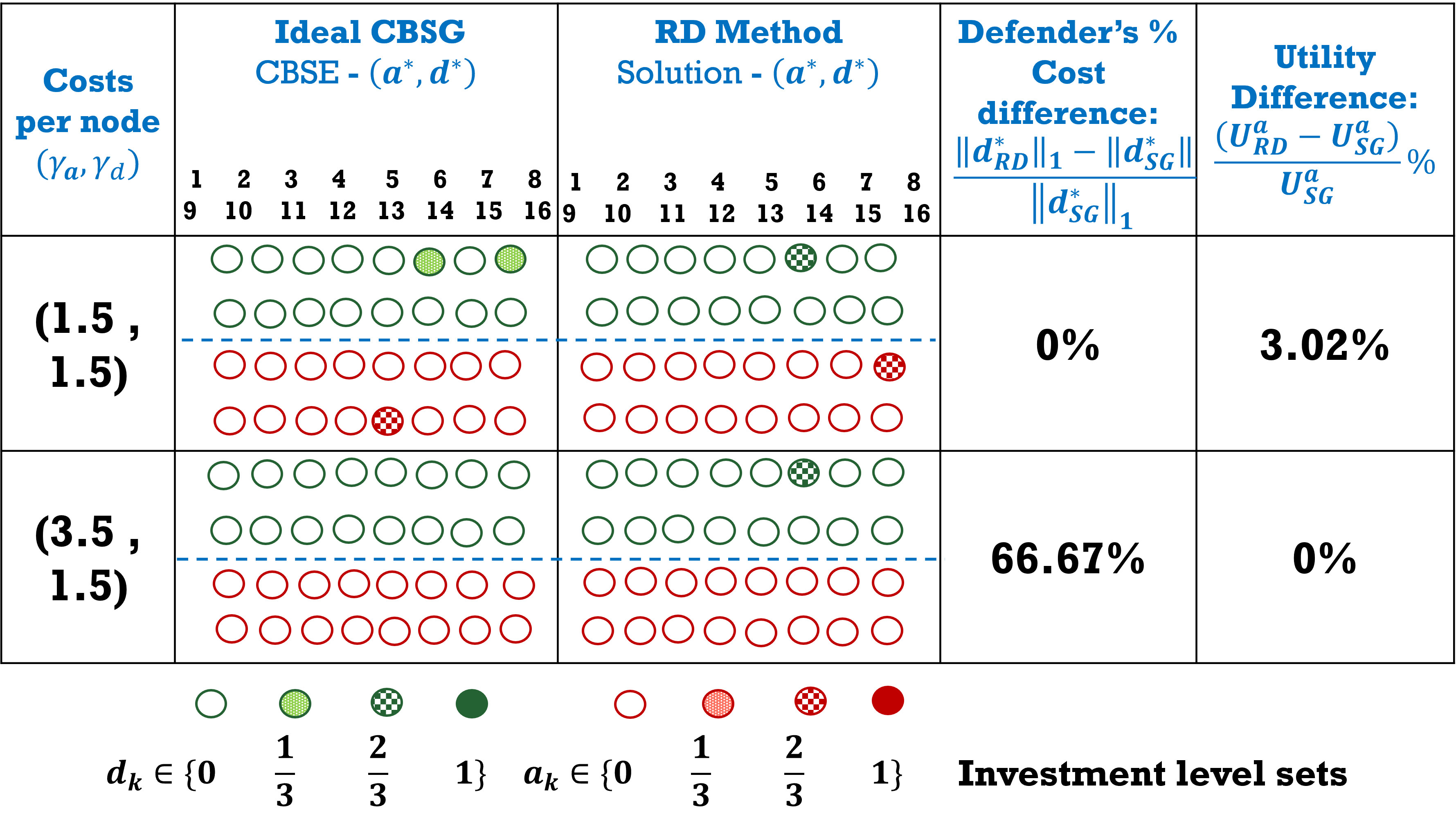}
       \caption*{(b)} 
      \end{minipage}
\caption{(a) Defender's fractional payoffs; comparison of the RD method (Section IV) and CBSG (Section III.C) vs. $\gamma_d$ for $\gamma_a = 1.5,2.5$;  (b) Players' strategies and the attacker's fractional payoff for two ($\gamma_a,\gamma_d)$ pairs; $L_a = L_d = 3$; the 68-bus system \cite{github_dataset}.}\label{fig17}
   \end{center}
\end{figure}

Next, we compare the performance of the proposed RD algorithm (see Section IV.B) to the ideal CBSG for the nominal-model system. Fig. \ref{fig17}a shows the payoffs of the defender vs. $\gamma_d$ for the ideal SG and the RD method for fixed values of cost of attack $\gamma_a$.  We observe the defender's utility  of the RD method does not exceed that of the ideal SG, consistent with Theorem 2,  but the former provides robust protection at small performance loss to the defender. In Fig.  \ref{fig17}b, we illustrate the strategies of the players  for these two methods and two cost pairs. These results as well as the mismatch trends \cite{PS_thesis} are consistent with those obtained for the IEEE 39-bus system in Fig. \ref{fig13}-\ref{fig14}. Furthermore, as in the IEEE 39-bus system, the boxplot statistics on the payoff difference for the RD method resemble those presented in Fig. \ref{fig16} for the nominal-payoff CBSG over the uncertain model set $\mathcal{M}$ \cite{PS_thesis}.

Finally, we discuss the computational complexity of the proposed methods for the IEEE 68-bus system. For  all algorithms of this subsection,  as for the IEEE 39-bus system, the computation of the loss vector corresponding to all sparsity patterns in (\ref{eq9}) is  dominated by the structured $\mathcal{H}_2$ optimization algorithm in \cite{Mihailo2013TAC}. This vector is computed in under $2500$ minutes$^1$. For Algorithm \ref{alg:ALG1}, the computation of the payoff matrix (\ref{eq13a}) is bounded by $(4^{16} \times 4^{16})$ for the SG designed for the IEEE 68-bus system model with $n=16$ generators and investment levels $L_a=L_d=3$ (Remark 3).  The traversal Algorithm \ref{alg:ALG1} (CBBI) runs out of memory space,  and the running time is exceeded when computing the payoff matrix for the cost pair $\gamma_a=\gamma_d = 0.5$, which corresponds to a very high computation load. Thus, when costs are small, the CBSG cannot be computed using Algorithm \ref{alg:ALG1}. On the other hand, the CB-BPEGA Algorithm \ref{alg:ALG2} with  the number of iterations $T=25$, population sizes $S_a =20$ and $S_d = 30$ for the attacker and defender, respectively, runs under $500$ seconds for any cost pair.  When the players' costs are moderate, for example, $\gamma_a=\gamma_d=1.5$, Algorithm \ref{alg:ALG1} (CBBI) runs in approximately $650$ seconds while Algorithm \ref{alg:ALG2} takes $490$ seconds to run. These results demonstrate efficiency of the evolutionary algorithms for large-scale systems.  The  complexity of the RD methods follows the trends discussed in Section V.B.5.

\section{Conclusion}
In this paper, we developed several investment methods for hardware or persistent DDoS attacks on NCSs that aim to degrade the control performance of the system. These methods allocate the resources of a system defender and malicious attacker strategically over the critical assets of an NCS under the players' limited budgets. First, a cost-based Stackelberg security-investment game between an attacker and a defender of an NCS  was investigated. The proposed SG allocate resources in a cost-efficient manner, and cost-based Stackelberg equilibria of the game reveal the ``important" nodes of the system, whose communication is critical for maintaining satisfactory control performance.  Moreover,  strategic investment approaches were proposed to address the model uncertainties of an NCS and defender's ignorance about the attacker's resources. To  reduce the complexity of the proposed  methods for large-scale systems, cost-based evolutionary methods were employed. Using an example of wide-area control of power systems applied to the IEEE 39-bus test model with model uncertainty arising from loads and generator parameters, we demonstrated that successful defense is feasible for both certain and uncertain scenarios unless the defender is much more resource-limited than the attacker. Furthermore, a case study of the IEEE 68-bus system was presented to validate the efficiency of the cost-based evolutionary methods for the proposed methods in large-scale networked control systems.

\appendix
\subsection{Proof of Theorem 1}
\noindent {\bf Theorem 1(a).} This proof follows from \cite{AnkurGameNotes}. Let $\Gamma^d$ denote the set of strategies $\bm{d}$ for the defender and $\Gamma^a$ the set of strategies $\bm{a}$ for the attacker. 
Since $\Gamma^d$ and $\Gamma^a$ are finite, the game is finite, i.e., the best response function $g(\bm{d})$ in Algorithm \ref{alg:ALG1} for any strategy $\bm{d}$ maps to a finite non-empty set (i.e. a maximum always exists for optimization (\ref{eq14}) when $\Gamma^a$ is finite). The defender's equilibrium strategy $\bm{d}^*$ exists since a maximum always exists for optimization (\ref{eq16}) when $\Gamma^d$ is finite. Similarly, the attacker's equilibrium strategy, given by $\bm{a}^*=g(\bm{d}^*)$, also exists. Thus, a strategy pair $(\bm{a}^*,\bm{d}^*)$ satisfying (\ref{eq14}), (\ref{eq16}) exists for the proposed~CBSG. \qed

\noindent {\bf Theorem 1(b).} Since the proposed CBSG is zero-sum, $U^a(\bm{a},\bm{d})=-U^d(\bm{a},\bm{d})$. The optimization problems (\ref{eq14}) and (\ref{eq16}) can be combined as a mini-max problem:
\begin{align}\label{eq28}
\bm{d}^* &= \mathop {\arg \max }\limits_{\bm{d}} U^d(g(\bm{d}),\bm{d})\nonumber \\
&=  - \arg \mathop {\max }\limits_{\bm{d}} U^a(g(\bm{d}),\bm{d})\\
&= \arg \mathop {\min }\limits_{\bm{d}} U^a(g(\bm{d}),\bm{d})\nonumber\\
&= \arg \mathop {\min }\limits_{\bm{d}} \mathop {\max }\limits_{\bm{a}} U^a(\bm{a},\bm{d}).\nonumber
\end{align}
Let the optimum of $\mathop {\min }\limits_{\bm{d}} \mathop {\max }\limits_{\bm{a}} U^a(\bm{a},\bm{d})$ be denoted as $U^a_o$.
Assume two different equilibrium strategy pairs $(\bm{a}^*_i, \bm{d}^*_i)$ and $(\bm{a}^*_j, \bm{d}^*_j)$, $i\neq j$, are solutions of (\ref{eq28}). Both strategy pairs must result in the same optimal payoff, i.e. $U^a(\bm{a}^*_i, \bm{d}^*_i)=U^a(\bm{a}^*_j, \bm{d}^*_j)=U^a_{o}$. Thus, the control performance loss associated with the payoffs, ${\bm{\Delta}_{\bm{s}_m}}$ in {\normalfont (\ref{eq12})} is the same for both~SEs. 

A similar argument can be made for the case when the game has $N$ SEs, where $N>2$. \qed

\noindent {\bf Theorem 1(c).} Assume $L_a$ and $L_d$ are fixed. When the attacker's cost per node is $\gamma_a$ and the defender's cost per node is $\gamma_d$, the defender's payoff at an SE is $U^d$. Let $\Gamma^a$ and $\Gamma^d$ denote the sets of strategies $\bm{a}$ and $\bm{d}$, respectively. Similarly, when the attacker's cost per node is $\gamma_a$, and the defender's cost per node increases to $\gamma_d'$, where $\gamma_d' > \gamma_d$, the defender's payoff at an SE is $U^{d'}$. In this case, the attacker's action space remains as $\Gamma^a$ while the defender's action space for $\gamma_d'$ is denoted $\Gamma^{d'}$. According to (\ref{eq11aa}), it is easy to show that $\Gamma^{d'} \subset \Gamma^d$ since $\gamma_d' > \gamma_d$.

Since CBSG is zero-sum, $U^a(\bm{a},\bm{d}) = - U^d(\bm{a},\bm{d})$. Eq. (\ref{eq14}) can be written as:
\begin{align}\label{eq29}
g(\bm{d}) &= \mathop {\arg \max }\limits_{\bm{a}} U^a(\bm{a},\bm{d})\nonumber \\
&=  \arg \mathop {\max }\limits_{\bm{a}} -U^d(\bm{a},\bm{d})\\
&= \arg \mathop {\min }\limits_{\bm{a}} U^d(\bm{a},\bm{d}).\nonumber
\end{align}
Combining (\ref{eq29}) with (\ref{eq16}), the optimization of the defender becomes:
\begin{align}\label{eq30}
\bm{d}^* &= \mathop {\arg \max }\limits_{\bm{d}} U^d(g(\bm{d}),\bm{d})\nonumber \\
&= \arg \mathop {\max }\limits_{\bm{d}} \mathop {\min }\limits_{\bm{a}} U^d(\bm{a},\bm{d}).
\end{align}
From (\ref{eq30}), given $\gamma_a$, the defender's payoff at an SE for $\gamma_d$ is computed as $U^d = \mathop {\max }\limits_{\bm{d}\in \Gamma^d} \mathop {\min }\limits_{\bm{a}\in \Gamma^a} U^d(\bm{a},\bm{d})$ while its payoff at an SE for $\gamma_d'$ is $U^{d'} = \mathop {\max }\limits_{\bm{d}'\in \Gamma^{d'}} \mathop {\min }\limits_{\bm{a}\in \Gamma^a} U^d(\bm{a},\bm{d})$. Since $\Gamma^{d'} \subset \Gamma^d$, $U^{d'} \leq U^d$. Thus, the defender's payoff is non-increasing with its cost per node $\gamma_d$ when $\gamma_a$ is fixed.

Similarly, we can make the argument that when the defender's cost per node $\gamma_d$ is fixed and the attacker's cost per node $\gamma_a$ increases, the attacker's payoff at an SE is non-increasing. \qed

\noindent {\bf Theorem 1(d)}
From (\ref{eq13a}), it is easy to show that the range of attacker's payoff at an SE is $[0, \Delta_{OL}]$, where $\Delta_{OL}$ is defined in (\ref{eq12l}). First, assume that $\epsilon < 1/n$, where $n$ is the total number of nodes in an NCS, and $\theta > L_d$. Then, from  (\ref{eq11aa}), when $\gamma_a<\epsilon$ and $\gamma_d>\theta$, the defender cannot protect any nodes while the attacker is able to attack all nodes at full effort, resulting in the  sparsity pattern $m = \bm{s}_{2^n -1}$ in (\ref{eq9}). Thus, feedback control is disabled and the attacker's utility at SE is $U^a(\bm{a}^*,\bm{d}^*) = \Delta_{OL}$.

Similarly, assume that $\alpha < 1/n$ and $\gamma_d < \alpha$. The defender can invest fully at each node following (\ref{eq11aa}) and  can protect against any attack in this case,  thus maintaining the optimal performance i.e., $J(\bm{K}^*_{\mathcal{H}_2})$ (dense $\mathcal{H}_2$ feedback control), resulting in the attacker's payoff $U^a(\bm{a}^*,\bm{d}^*) = 0$ at an SE. \qed

\noindent {\bf Theorem 1(e).} Assume $\gamma_a$ and $\gamma_d$ are fixed. Let $U^d$ denote the defender's payoff at an SE when the attacker's number of investment levels is $L_a+1$ and the defender's number of investment levels is $L_d+1$.  Similarly, when the attacker's number of investment levels is $L_a+1$ and the defender's number of investment levels is $L_d'+1$, let $U^{d'}$ denote the defender's payoff at an SE.
 
When the defender's number of investment levels is $L_d+1$, the defender's actions are represented by $\bm{d}=(d_1,...,d_k,...,d_n)$, where $d_k\in \bm{D} = \{0,\frac{1}{L_d},\frac{2}{L_d},\cdots,1\}$ denotes the defender's investment level into node $k$. Let $\Gamma^d$ denote the set of strategies $\bm{d}$ for the defender for this case.
 
When the defender's number of investment levels is increased to $L_d'=\eta(L_d)$, where $\eta$ is a positive integer, the set of possible values of $d_k \in \bm{D}'= \{ 0,\frac{1}{L_d'},\frac{2}{L_d'},\cdots,1\}$. Let $\Gamma^{d'}$ denote the set of strategies $\bm{d}'$ for the defender for this case. Since $L_d'=\eta(L_d)$, it is easy to show that $\bm{D} \subset \bm{D}'$. Thus, $\Gamma^d \subset \Gamma^{d'}$.
 
 When computing the SEs, the defender searches all possible actions in its action space $\Gamma^d$ (for $L_d$) or $\Gamma^{d'}$ (for $L_d'$). Since $\Gamma^d \subset \Gamma^{d'}$, $U^d \leq U^{d'}$. Thus, the defender's payoff does not decrease when $L_d$ is increased to $L_d'=\eta(L_d)$, where $\eta$ is a positive integer given the costs and $L_a+1$ are fixed.
 
Similar argument can be made for the case when $L_d+1$ is fixed for the defender and $L_a$ is increased to a larger number $L_a'=\eta(L_a)$ for the attacker. \qed

\subsection{Proof of Proposition 1}
First, since the action spaces of both players are finite, from Theorem 1(a), an SE exists in the proposed game. Moreover, according to Theorem 1(b), the players' payoffs are the same at any CBSE of the proposed zero-sum CBSG.

From Theorem 4.1 in \cite{GA_proof}, due to the crossover and mutation process, any feasible solution appears at least once in some generation with probability 1 when the set of feasible solutions is finite and the maximum number of iterations $T$ tends to infinity. Thus, since the action spaces of both players in the proposed game are finite, any feasible strategy pair $(\bm{a},\bm{d})$ appears at least once in some generation with probability 1 when $T \rightarrow \infty$.

 Denote the set of all SEs as $\mathcal{E} = \left\{(\bm{a}^*_i,\bm{d}^*_i) : i = 1,...,S\right\}$. Suppose $(\bm{a}_o^*,\bm{d}_o^*)$ is a CBSE for the proposed zero-sum CBSG, i.e., this strategy pair is an SE with the lowest players' costs satisfying (\ref{eq11aa}). Once the least-cost strategy pair $(\bm{a}_o^*,\bm{d}_o^*)$ appears in some generation, its fitness value (\ref{eq21fitd}) (or (\ref{eq22fita})) for the defender (or attacker) must be the highest among the members of that generation since all SEs have the same payoff for the players in the proposed zero-sum SG. Moreover, since $\bm{d}^*_o$ (or $\bm{a}^*_o$) has the smallest cost among all $\bm{d}^*_i$ (or $\bm{a}^*_i$), it must have the highest ranking among the members of that generation due to the sorting and reordering property in Step 6 of Algorithm \ref{alg:ALG2}. Thus, $\bm{d}^*_o$ (or $\bm{a}^*_o$) will be included in the next generation. 
 
 By the argument above, the CBSE pair (least-cost SE) $(\bm{a}_o^*,\bm{d}_o^*)$ will occur in all subsequent generations of CB-BPEGA. Thus, it will be included in the final generations $POP_a^T$ and $POP_d^T$ as $T \rightarrow \infty$.
 
 Finally according to Step 7 of Algorithm \ref{alg:ALG2}, the CBBI algorithm chooses the least-cost SE, i.e. CBSE $(\bm{a}_o^*,\bm{d}_o^*)$ from the set of SEs in the final generations. \qed

\subsection{Proof of Theorem 2}
Let $\Gamma^a$ and $\Gamma^d$ denote the action sets for the attacker and defender, respectively, for costs of attack and defense per node as $\gamma_a$ and $\gamma_d$ {\normalfont (\ref{eq11aa})}. Let $(\bm{a}^*,\bm{d}^*)$ denote an CBSE of the ideal CBSG  {\normalfont (\ref{eq13a})-(\ref{eq13b})} while $(\bm{a}^*_{RD},\bm{d}^*_{RD})$ denotes the solution of the RD method {\normalfont (\ref{eqactualpayoffrdsg})}.
 
\noindent Then the following conditions hold
 \begin{align}
     U^a(\bm{a}^*_o,\bm{d}^*_o) \geq U^a(\bm{a},\bm{d}^*_o), \ \forall \bm{a}\in \Gamma^a \\
      U^d(\bm{a}^*_o,\bm{d}^*_o) \geq U^d(\bm{a}^*_o,\bm{d}), \ \forall \bm{d}\in \Gamma^d \label{36}
 \end{align}
 Since $\bm{d}^*_{RD} \in \Gamma^d$, from {\normalfont (\ref{36})} 
 \begin{equation}
     U^d(\bm{a}^*_O,\bm{d}^*_o) \geq U^d(\bm{a}^*_o,\bm{d}^*_{RD})
 \end{equation}
 Since the CBSG  {\normalfont (\ref{eq13a})-(\ref{eq13b})} is zero-sum,
 \begin{equation}
     U^a(\bm{a}^*_o,\bm{d}^*_o) \leq U^a(\bm{a}^*_o,\bm{d}^*_{RD}) \label{38}
 \end{equation}
 Furthermore, the attacker's payoff at the solution of  RD method {\normalfont (\ref{eqactualpayoffrdsg})} satisfies:
 \begin{align}
 U^a(\bm{a}^*_{RD},\bm{d}^*_{RD}) \geq U^a(\bm{a},\bm{d}^*_{RD}) \ \forall \bm{a} \in \Gamma^a \label{40}
 \end{align}
 Thus,
 \begin{equation}
     U^a(\bm{a}^*_{RD},\bm{d}^*_{RD}) \geq U^a(\bm{a}^*_o,\bm{d}^*_{RD}) \label{41}
 \end{equation}
 From {\normalfont (\ref{38})}  and {\normalfont (\ref{41})} ,
 \begin{equation}
     U^a(\bm{a}^*_o,\bm{d}^*_o) \leq U^a(\bm{a}^*_{RD},\bm{d}^*_{RD})
 \end{equation}
 and
 \begin{equation}
     U^d(\bm{a}^*_o,\bm{d}^*_o) \geq  U^d(\bm{a}^*_{RD},\bm{d}^*_{RD})
 \end{equation} \qed

\subsection{Model Parameters}
Matrices for the linearized state-space model in \ref{eq1} are provided in \cite{github_dataset} for both case studies, where $\bm{A} \in \mathbb{R}^{75 \times 75}$,  $\bm{B} \in \mathbb{R}^{75 \times 9}$ and $\bm{K} \in \mathbb{R}^{9 \times 75}$ for the IEEE 39-bus system and $\bm{A}\in \mathbb{R}^{160 \times 160}$, $\bm{B}\in \mathbb{R}^{160 \times 16}$ and $\bm{K} \in \mathbb{R}^{16 \times 160}$ for the IEEE 68-bus system. The nonlinear data file containing bus, line, load, and generator information is also attached in \cite{github_dataset}. Finally, the data files containing state-space matrices for linearized uncertain models generated for each of the system can also be found in \cite{github_dataset}.

\balance

\ifCLASSOPTIONcaptionsoff
  \newpage
\fi

 \bibliographystyle{IEEEtran}
\bibliography{references}

\end{document}